\documentclass[11pt]{amsart}
\usepackage{amssymb}
\usepackage{amsmath}
\usepackage{amsrefs}
\usepackage[colorlinks=true]{hyperref}
\usepackage{cleveref}
\usepackage{enumitem}
\usepackage{amscd}
\usepackage{amsthm}
\usepackage{mathtools}
\usepackage{graphicx}
\usepackage{tikz}
\usetikzlibrary{decorations.pathreplacing,decorations.markings,hobby,knots,shapes.geometric,calc}

\usepackage{mathrsfs}
\usepackage{pinlabel}
\usepackage[margin=2.93cm]{geometry}

\usepackage{xcolor}
\hypersetup{
    colorlinks,
    linkcolor={red!50!black},
    citecolor={green!50!black},
    urlcolor={blue!80!black}
}
\allowdisplaybreaks

\newtheorem{theorem}{Theorem}[section]
\newtheorem{lemma}[theorem]{Lemma}
\newtheorem{proposition}[theorem]{Proposition}
\newtheorem{corollary}[theorem]{Corollary}

\newtheorem{conjecture}[theorem]{Conjecture}
\newtheorem{question}[theorem]{Question}
\newtheorem{assumption}[theorem]{Assumption}

\newtheorem*{assumption-no-number}{Assumption}
\newtheorem*{corollary*}{Corollary}

\theoremstyle{definition}
\newtheorem{definition}[theorem]{Definition}
\newtheorem{example}[theorem]{Example}

\theoremstyle{remark}
\newtheorem{remark}[theorem]{Remark}

\numberwithin{equation}{section}

\newcommand{\R}{\mathbb{R}}

\newcommand{\N}{\mathbb{N}}
\newcommand{\Z}{\mathbb{Z}}

\renewcommand{\L}{\mathcal{L}}
\newcommand{\pt}{\mathrm{pt}}

\newcommand{\lt}{\left}
\newcommand{\rt}{\right}
\newcommand{\tms}{\times}

\newcommand{\rmk}{\begin{remark}}
\newcommand{\ermk}{\end{remark}}
\newcommand{\cor}{\begin{corollary}}
\newcommand{\ecor}{\end{corollary}}
\newcommand{\eq}{\begin{equation}}
\newcommand{\eeq}{\end{equation}}
\newcommand{\eqs}{\begin{equation*}}
\newcommand{\eeqs}{\end{equation*}}
\newcommand{\prop}{\begin{proposition}}
\newcommand{\eprop}{\end{proposition}}
\newcommand{\thm}{\begin{theorem}}
\newcommand{\ethm}{\end{theorem}}
\newcommand{\conj}{\begin{conjecture}}
\newcommand{\econj}{\end{conjecture}}
\newcommand{\lem}{\begin{lemma}}
\newcommand{\elem}{\end{lemma}}
\newcommand{\defi}{\begin{definition}}
\newcommand{\edefi}{\end{definition}}
\newcommand{\ex}{\begin{example}}
\newcommand{\eex}{\end{example}}
\newcommand{\alis}{\begin{align*}}
\newcommand{\ealis}{\end{align*}}
\newcommand{\pf}{\begin{proof}}
\newcommand{\epf}{\end{proof}}
\newcommand{\ali}{\begin{align}}
\newcommand{\eali}{\end{align}}
\newcommand{\qus}{\begin{question}}
\newcommand{\equs}{\end{question}}

\newcommand{\mc}{\mathcal}
\renewcommand{\bf}{\textbf}

\newcommand{\C}{\mathbb{C}}
\newcommand{\TT}{\mathbb{T}}

\newcommand{\sub}{\subset}

\newcommand{\ov}{\overline}

\newcommand{\bb}{\mathbb}

\newcommand{\op}{\operatorname}

\renewcommand{\a}{\alpha}
\renewcommand{\b}{\beta}
\renewcommand{\d}{\partial}
\newcommand{\e}{\epsilon}

\newcommand{\g}{\gamma}

\newcommand{\s}{\sigma}
\renewcommand{\t}{\theta}

\renewcommand{\l}{\lambda}
\renewcommand{\o}{\omega}
\newcommand{\z}{\zeta}

\newcommand{\G}{\Gamma}
\renewcommand{\L}{\Lambda}
\renewcommand{\O}{\Omega}

\renewcommand{\S}{\Sigma}

\newcommand{\Q}{\mathbb{Q}}

\renewcommand{\ov}{\overline}

\begin{document}
\title[Linking of Tori and Embedding Obstructions]{Linking of Lagrangian Tori and Embedding Obstructions in Symplectic 4-manifolds}
\author[Laurent C\^ot\'e]{Laurent C\^ot\'e}
\thanks {LC was supported by a Stanford University Benchmark Graduate Fellowship.}
\address{Department of Mathematics, Stanford University, 450 Serra Mall, Stanford, CA 94305}
\email{lcote@stanford.edu}

\author[Georgios Dimitroglou Rizell]{Georgios Dimitroglou Rizell}
\thanks {GDR is supported by the grant KAW 2016.0198 from the Knut and Alice Wallenberg Foundation.}
\address {Department of Mathematics, Uppsala University, Box 480, SE-751 06, Uppsala, Sweden}
\email {georgios.dimitroglou@math.uu.se}

\begin{abstract} 
We classify weakly exact, rational Lagrangian tori in $T^* \bb{T}^2- 0_{\bb{T}^2}$ up to Hamiltonian isotopy. This result is related to the classification theory of closed $1$-forms on $\bb{T}^n$ and also has applications to symplectic topology. As a first corollary, we strengthen a result due independently to Eliashberg--Polterovich and to Giroux describing Lagrangian tori in $T^* \bb{T}^2-0_{\bb{T}^2}$ which are homologous to the zero section. As a second corollary, we exhibit pairs of disjoint totally real tori $K_1, K_2 \sub T^*\bb{T}^2$, each of which is isotopic through totally real tori to the zero section, but such that the union $K_1 \cup K_2$ is not even smoothly isotopic to a Lagrangian.   

In the second part of the paper, we study linking of Lagrangian tori in $(\R^4, \o)$ and in rational symplectic $4$-manifolds. We prove that the linking properties of such tori are determined by purely algebro-topological data, which can often be deduced from enumerative disk counts in the monotone case. We also use this result to describe certain Lagrangian embedding obstructions. \end{abstract} 
\maketitle

\section{Introduction}


\subsection{Lagrangian tori in $T^*\bb{T}^2$ in the complement of the zero section}

Let $\bb{T}^2= \R^2/ \Z^2$ be the $2$-torus. We consider its cotangent bundle $(T^* \bb{T}^2, d\l_{\op{can}})= (\R^2/\Z^2 \tms \R^2, y_1 dx_1 + y_2 dx_2)$. 

The first theorem of this paper provides a classification of weakly exact, rational Lagrangian tori in $T^*\bb{T}^2-0_{\bb{T}^2}$ up to Hamiltonian isotopy. Recall that a Lagrangian $L$ in a symplectic manifold $(M, \o)$ is said to be \emph{weakly exact} if $\o(\a) =0$ for all $\a \in \pi_2(M, L)$. It is said to be \emph{rational} if $\l|_L \in H^1(L; \R)$ is contained in the subset $c \cdot H^1(L; \Q) \sub H^1(L; \R)$, for some $c \in \R$. 

\thm \label{theorem:weakly-exact-classification} Let $L \sub T^*\bb{T}^2-0_{\bb{T}^2}$ be a weakly exact, rational Lagrangian torus. Then $L$ is Hamiltonian isotopic in $T^*\bb{T}^2-0_{\bb{T}^2}$ to a product torus of the form $\R^2/\Z^2 \tms \{\mathbf{p}\}$, for some $\mathbf{p} \neq (0, 0)$. 
\ethm 

\Cref{theorem:weakly-exact-classification} can be motivated from two perspectives. On the one hand, it has several applications to symplectic topology. On the other hand, it is related to the problem of classifying closed $1$-forms on tori, which was studied by various authors including Laudenbach \cite{laudenbach}, Laudenbach--Blank \cite{laudenbach-blank} and Sikorav \cite{sikorav2}. We discuss both of these perspectives in more detail for the remainder of this section.

Our first application of \Cref{theorem:weakly-exact-classification} to symplectic topology is following corollary, which is proved in \Cref{subsection:some-background}.

\cor  \label{corollary:eli-pol-giroux} Let $L \sub T^*\bb{T}^2 = \bb{T}^2 \tms \R^2$ be a Lagrangian torus which is disjoint from the zero section and homologous to it. Then $L$ is isotopic through Lagrangian tori disjoint from the zero section to $\bb{T}^2 \tms \{\mathbf{p} \}$, for some $\mathbf{p} \neq (0,0)$.
\ecor  

\Cref{corollary:eli-pol-giroux} strengthens a result due independently to Eliashberg--Polterovich \cite{eli-pol} and to Giroux \cite{giroux}, who proved under the same hypotheses that the projection $L \to \R^2 - \{0\}$ is homotopic to a point. 

\subsubsection{Obstructions to Lagrangian embeddings in $T^* \bb{T}^2$} \label{subsubsection:embeddings-1}

\Cref{theorem:weakly-exact-classification} can be shown to imply the following technical statement, which is proved in \Cref{subsection:some-background}.

\cor \label{corollary:pi-1-cotangent}
	Let $L_0, L_1 \sub T^*\bb{T}^2$ be disjoint Lagrangian tori. Consider the subgroup $\G_j \sub \pi_1(T^*\bb{T}^2- L_j)$ defined by setting $\G_j= \op{Im}( \pi_1(\bb{T}^2 \tms \{ \mathbf{p} \}) \to \pi_1(T^*\bb{T}^2- L_j))$ for $\|\mathbf{p} \|$ sufficiently large. Suppose that the inclusions $L_i \hookrightarrow T^*\bb{T}^2$ are homotopy equivalences.  Then the image of the natural map $\pi_1(L_i) \to \pi_1(T^*\bb{T}^2- L_j)$ coincides with $\G_j$ for $i \neq j$. 
\ecor

\begin{figure}[htp]
\centering
\vspace{3mm}
\labellist
\pinlabel $\color{red}\ell_2$ at 37 5
\pinlabel $\color{blue}\ell_1$ at -1 5
\endlabellist
\includegraphics[scale=4]{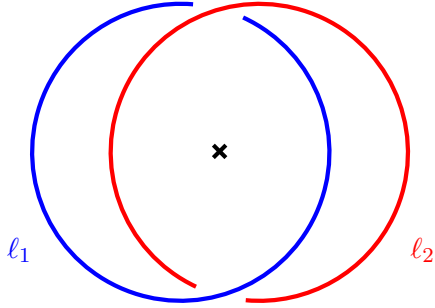}
\caption{A homotopically essential Hopf link inside $\C^*\times \R$; the picture shows the projection onto the $\C^*$ component.}
\label{picture:hopf}
\end{figure}

\Cref{corollary:pi-1-cotangent} gives new obstructions to Lagrangian embeddings in $T^* \bb{T}^2$. Fixing the standard integrable complex structure on $T^*\bb{T}^2$, let $K_1, K_1 \sub \C^* \tms \R \tms S^1 = T^* \bb{T}^2$ be the totally real tori obtained by ``spinning" the Hopf link $\ell= \ell_1 \cup \ell_2$ drawn in \Cref{picture:hopf}; i.e. by taking the product with $S^1$.  

It is clear that both $K_0$ and $K_1$ are isotopic through totally real tori to a Lagrangian torus (indeed, both tori are isotopic to the zero section). However, we claim that the disconnected manifold $K_0 \cup K_1$ is not smoothly isotopic to a Lagrangian embedding of $\bb{T}^2 \sqcup \bb{T}^2$. Indeed, for $p \in S^1$, the class $[\ell_i \tms p]$ is in the image of the natural map $\pi_1(K_i) \to \pi_1(T^*\bb{T}^2- K_j)$. However, this class is not in the image of $\pi_1(\bb{T}^2 \tms \{\mathbf{p} \}) \to \pi_1(T^*\bb{T}^2- K_j)$ for $\| \mathbf{p} \|$ large. It follows by \Cref{corollary:pi-1-cotangent} that $K_0 \cup K_1$ is not isotopic to a Lagrangian embedding of $\bb{T}^2 \sqcup \bb{T}^2$. 

We summarize this as a corollary. 

\cor \label{corollary:tori-linked-cotangent}
There exist pairs of disjoint, totally real tori $K_1, K_2 \sub T^*\bb{T}^2$ with the property that $K_1, K_2$ are each isotopic to the zero section through totally real tori, but such that $K_1 \cup K_2$ is not smoothly isotopic to a Lagrangian submanifold. 
\ecor

This phenomenon is somewhat analogous to recent work of Mak--Smith \cite{mak-smith}, who exhibited a pair of displaceable, disjoint Lagrangian tori $L_0, L_1 \sub (S^2 \tms S^2, A\o \oplus a\o)$, with $0<a<A$, such that $L_0 \cup L_1$ is not displaceable. We will describe other examples in $\R^4, \bb{CP}^2$ and $S^2 \tms S^2$ in \Cref{subsection:obstructions-linked-r4}. 

Note that \Cref{corollary:tori-linked-cotangent} and the example preceding it can easily be generalized to more complicated links, with an arbitrary number of components.

\subsubsection{Classification theory of closed $1$-forms}
Following Laudenbach \cite{laudenbach}, let us say that two nowhere vanishing closed $1$-forms on $\bb{T}^n$ are \emph{equivalent} if they can be joined by a path of cohomologous nowhere vanishing closed $1$-forms.  One then naturally seeks to classify cohomologous nowhere vanishing closed $1$-forms up to equivalence. 

When $n=2$, it can be shown using the Poincar\'{e}--Bendixon theorem that all cohomologous nonvanishing closed $1$-forms are equivalent. When $n=3$, this is a result due to Laudenbach--Blank \cite{laudenbach-blank}. In contrast, for $n \geq 6$, Laudenbach showed \cite{laudenbach} that there are infinitely many non-equivalent nonvanishing $1$-forms on $\bb{T}^n$ representing the same rational cohomology class. 

This dichotomy can be explained from the perspective of algebraic $K$-theory. Indeed, Laudenbach exhibited an obstruction to rational cohomologous $1$-forms on $\bb{T}^n$ being equivalent which lives in the group $P(\bb{T}^{n-1})$ of pseudo-isotopy classes of $\bb{T}^{n-1}$. This group is trivial for $n=2,3$. In contrast, for $n \geq 6$, a celebrated theorem of Hatcher and Wagoner \cite{hatcher-wagoner} implies that $P(\bb{T}^{n-1})$ surjects onto $\op{Wh}_2(\bb{T}^{n-1})$. In fact, $P(\bb{T}^{n-1})$ is isomorphic as a group to the algebra of sets of $\Z^n- \{0\}$. 

Observe that a nowhere vanishing closed $1$-form on $\bb{T}^n$ is the same thing as a \emph{graphical} Lagrangian $n$-torus in $T^* \bb{T}^n- 0_{\bb{T}^n}$. The classification of nowhere vanishing closed $1$-forms on $\bb{T}^n$ up to equivalence can thus be understood as the classification of \emph{graphical} Lagrangian tori, up to Hamiltonian isotopy though graphical Lagrangian tori. From this perspective, \Cref{theorem:weakly-exact-classification} generalizes the classification of graphical Lagrangian tori in $T^*\bb{T}^2- 0_{\bb{T}^2}$ to the non-graphical case, under a rationality assumption.  

In light of Laudenbach's results mentioned above, \Cref{theorem:weakly-exact-classification} appears to reflect the fact that the pseudo-isotopy group of the circle is trivial. In particular, one does not expect \Cref{theorem:weakly-exact-classification} to generalize to dimensions $\geq 6$. One potential approach to proving this would be to develop $K$-theoretic invariants in Floer theory for possibly non-exact Lagrangians (see \cite{sullivan} for related work in the exact setting).

\subsection{Linking of Lagrangian tori in rational symplectic manifolds}

In the second part of this paper, we study linking of Lagrangian tori in $(\R^4, \o)$ and in rational symplectic $4$-manifolds. Following \cite{lalonde-mcduff}, we say that a symplectic $4$-manifold is \emph{rational} if it can be obtained from $\bb{CP}^2$ by a sequence of blowups and blowdowns; see\ \Cref{subsection:unlinking-rational}. 

Recall that a Lagrangian $L \sub (M, \o)$ is said to be \emph{monotone} with monotonicity constant $\kappa>0$ if $\o = \kappa \mu \in \pi_2(M, L)$, where $\mu$ is the Maslov class. It will be convenient to make the following definition.

\defi
	Given a topological space $X$ and a subspace $\S \sub X$, we say that $\S$ is \emph{contractible} if the inclusion $\S \hookrightarrow X$ is homotopic to a constant map. 
\edefi

It was shown in \cite[Thm.\ A]{dgi} that all Lagrangian tori in $\R^4, S^2 \tms S^2, \bb{CP}^2$ (with their standard symplectic structures) are Lagrangian isotopic. In particular, since any Darboux chart contains a Lagrangian torus, all such tori are contractible. We prove the following extension to all rational symplectic $4$-manifolds (see \Cref{subsection:unlinking-rational}).

\prop  \label{proposition:all-tori-isotopic}
	All Lagrangian tori in a rational symplectic $4$-manifold are smoothly isotopic. In particular, they are all contractible.
\eprop

Given a symplectic manifold $(M,\o)$ and a collection $L_1,\dots,L_n \sub M$ of disjoint, contractible Lagrangian submanifolds, we say that the $L_i$ are \emph{unlinked} if they can be isotoped into disjoint balls without intersecting. More precisely, let us consider the following definition (cf.\cite[Def.\ 1.1]{cote}): 

\defi
Let $L_1,\dots,L_n$ be disjoint, contractible Lagrangians in a symplectic manifold $(X, \o)$. We say that the $L_i$ are \emph{smoothly unlinked} if there exists a collection of disjoint embedded closed balls $B_1,\dots,B_n \sub X$ and $1$-parameter family of smooth embeddings $\Phi: [0,1] \tms (\bigsqcup_{i=1}^n L_i ) \to X$ such that, for $i=1,\dots,n$, the composition $$L_i \xhookrightarrow{L_i \mapsto \{t\} \tms L_i} [0,1] \tms \lt(\bigsqcup_{i=1}^n L_i\rt) \xrightarrow{\Phi} X$$ maps $L_i$ to $L_i \sub X$ diffeomorphically if $t=0$ and maps $L_i$ into $B_i$ if $t=1$.\footnote{We are slightly abusing notation by viewing $L_i$ both as an abstract manifold and as an embedded submanifold of $X$.}

We say that the $L_i$ are \emph{Lagrangian unlinked} (resp. \emph{Hamiltonian unlinked}) if one can choose $\Phi$ so that $\Phi(t,-)$ is a Lagrangian embedding for all $t \in [0,1]$ (resp. if $\Phi$ is induced by a global Hamiltonian isotopy of $(X,\o)$). 
\edefi 

The following theorem shows that Lagrangian tori in $\R^4$ are Lagrangian unlinked if and only if an obviously necessary algebro-topological condition is satisfied. 

\thm \label{theorem:lagrangian-r4-introduction}
	Let $L_1,\dots,L_n$ be disjoint Lagrangian tori in the symplectic vector space $(\R^4, \o= dx_1\wedge dy^1 + dx_2 \wedge dy_2)$. If the inclusion 
	\eq \label{equation:inclusion} \iota_i: L_i \to \R^4- \cup_{j \neq i} L_j\eeq induces the trivial map on fundamental groups for $i=1,2\dots,n$, then the $L_i$ are Lagrangian unlinked.  
	
	If the $L_i$ are monotone and $\kappa_i \leq \kappa_j$ for $i \leq j$, where $\kappa_i$ is the monotonicity constant of $L_i$, then the conclusion also holds if we only assume that $\tilde{\iota}_i: L_i \to \R^4- \bigcup_{j >i} L_j$ induces the trivial map on fundamental groups.
\ethm

We also prove a weaker version of \Cref{theorem:lagrangian-r4-introduction} for rational symplectic $4$-manifolds.

\thm \label{theorem:lagrangian-rational-introduction}
	Let $(X, \o)$ be a rational symplectic $4$-manifold and let $L_1,\dots,L_n \sub (X, \o)$ be disjoint Lagrangian tori. Suppose that the inclusion 
	\eq \label{equation:inclusion-2} \iota_i: L_i \to X- \cup_{j \neq i} L_j \eeq induces the trivial map on fundamental groups for $i=1,\dots,n$. Then the $L_i$ are smoothly unlinked. If moreover $(X, \o)$ is minimal, then the $L_i$ are Lagrangian unlinked. 
%
\ethm

The proofs of Theorems \ref{theorem:lagrangian-r4-introduction} and \ref{theorem:lagrangian-rational-introduction} are respectively provided in Sections \ref{subsection:unlinking-r4} and \ref{subsection:unlinking-rational}.

Observe that the algebro-topological conditions \eqref{equation:inclusion} and \eqref{equation:inclusion-2} are clearly necessary for conclusions of the above theorems to hold. However, it is not at all obvious that these conditions are sufficient. Indeed, there could a priori be additional obstructions to unlinking the tori coming from smooth topology or symplectic topology. 

We remark that it is not difficult to construct configurations of Lagrangian tori which are linked, both in $\R^4$ and in rational symplectic manifolds (see for instance \cite[Ex.\ 4.10]{cote}). 

Building on work of the first author \cite{cote}, one can show that \eqref{equation:inclusion} or \eqref{equation:inclusion-2} are automatically satisfied for monotone Lagrangian tori with ``enough" nonvanishing counts of Maslov $2$ disks.  

To make this precise, let $n(L, \beta)$ be the mod-$2$ count of $J$-holomorphic disks passing through a generic point of $L$, for a generic $J$, where $L \sub (M, \o)$ is a monotone Lagrangian and $\beta \in \pi_2(M, L)$ has Maslov number $2$. We then have the following proposition.

\prop  \label{prop:monotone-intro} Let $L_1,\dots,L_n \sub (M, \o)$ be disjoint monotone Lagrangian tori, where $L_i$ has monotonicity constant $\kappa_i$. Suppose that there exists a collection of classes $\mc{C}_i = \{ \beta_1,\dots,\beta_{n_i} \} \sub \pi_2(M, L_i)$ such that $n(L_i, \beta_k) \neq 0$ for all $1 \leq k \leq n_i$, and such that the image of $\mc{C}_i$ under the composition 
\eq  \pi_2(M, L_i) \to H_2(M, L_i) \to H_1(L_i) = \pi_1(L_i) \eeq generates $\pi_1(L_i)$ as an abelian group. Then $\iota_i: L_i \to  M- \bigcup_{j>i} L_j$ induces the trivial map on fundamental groups. 
\eprop

If we consider for example a collection $L_1,\dots,L_n \sub (\R^4, \o)$ where the $L_i$ are Clifford tori, then \Cref{prop:monotone-intro} combined with \Cref{theorem:lagrangian-r4-introduction} implies that the $L_i$ are Lagrangian unlinked. This strengthens \cite[Thm.\ A]{cote} of the first author, who showed that the $L_i$ are smoothly unlinked. By combining \Cref{prop:monotone-intro} with \Cref{theorem:lagrangian-rational-introduction}, one obtains similar conclusions for the ``exotic" tori in $\bb{CP}^2$ or Del Pezzo surfaces such as those constructed by Vianna \cite{vianna, vianna2}. However, it's unclear whether one can find collections of such tori which do not pairwise intersect, so the result may be vacuous in this case. 

\subsubsection{Obstructions to Lagrangian embeddings} \label{subsection:obstructions-linked-r4}

Theorems \ref{theorem:lagrangian-r4-introduction} and \ref{theorem:lagrangian-rational-introduction} can be used to obstruct certain Lagrangian embeddings of $\bb{T}^2 \sqcup \bb{T}^2$ into $\R^4, S^2 \tms S^2$ and $\bb{CP}^2$. 

Let $\ell= \ell_1 \cup \ell_2 \sub \bb{H}= \{(x,y, z) \mid z >0 \}$ be a two-component link satisfying the following properties:
\begin{itemize}
\item[(i)] $\pi_1(\bb{H} - \ell ) $ is not equal to $\Z * \Z$, i.e. the fundamental group of the link complement is distinct from that of the unlink. 
\item[(ii)] $\ell_i \sub \bb{H}- \ell_j$ is nullhomotopic for $i\neq j$.
\item[(iii)] $\ell_i$ is smoothly isotopic as an embedded submanifold to $\{(x,y,z) \mid x^2+y^2=1, z=1\}$, via a family of embeddings which are never tangent to the lines $\{x=\op{constant}, y=\op{constant} \}$.
\end{itemize}

\begin{figure}[htp]
\centering
\vspace{3mm}
\labellist
\pinlabel $\color{blue}\ell_1$ at -4 23
\pinlabel $\color{red}\ell_2$ at 22 23
\endlabellist
\includegraphics[scale=3]{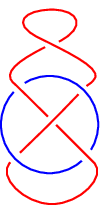}
\caption{A Whitehead link inside $\bb{H}= \{(x,y, z) \mid z >0 \}$; the picture shows the projection onto the $(x,y)$-component}
\label{picture:whitehead}
\end{figure}

An example of such a link is the Whitehead link drawn in \Cref{picture:whitehead}. The fundamental group of the Whitehead link complement is computed in \cite[Lem.\ 9.4]{cooper-long} and can be shown to be different from that of the unlink (this follows, for instance, from the fact that the $\op{SL}(2,\C)$ character variety of the Whitehead link complement has a $2$-dimensional component \cite[Appendix C]{tillmann} while the character variety of the free group is isomorphic to $\C^3$ \cite[Sec.\ 2.3]{abou-man}). 

To verify (ii), note that one component of the Whitehead link is clearly nullhomotopic in the complement of the other. However, one can also show that the link projection \Cref{picture:whitehead} is symmetric in the two link components, i.e. there is an ambient isotopy which interchanges $\ell_1$ and $\ell_2$. 

Finally, (iii) is obvious from \Cref{picture:whitehead}. 

For $i=1,2$, let $K_i \sub \bb{H} \tms S^1$ be obtained by ``spinning" the link component $\ell_i$. We view the $K_i$ as tori in $\R^4$ via the embedding $ \bb{H} \tms S^1 \to \R^4$ taking $(x,y, z, \t) \mapsto (x, y, z \cos \t, z \sin \t)$. It is straightforward to check using the van Kampen theorem and (i) that $\pi_1(\R^4-(K_1 \cup K_2))= \pi_1(\bb{H}- (\ell_1 \cup \ell_2)) \neq \Z* \Z$. 

Observe that the $K_i$ are totally real with respect to the standard complex structure on $\R^4$. Moreover, it follows from (iii) that each $K_i$ is isotopic through totally real tori to the standard Clifford torus $\{(z_1, z_2) \in \C^2 \mid |z_1|= |z_2|= 1 \}$.

We claim that $K_1 \cup K_2$ is not smoothly isotopic to a Lagrangian embedding of $\bb{T}^2 \sqcup \bb{T}^2$ into $\R^4$. Indeed, it follows easily from (ii) that the inclusion $K_i \to \R^4-K_j$ induces the trivial map on fundamental groups for $i \neq j$. Hence \Cref{theorem:lagrangian-r4-introduction} implies that $K_1, K_2$ are Lagrangian unlinked. Hence $\pi_1(\R^4- (K_1 \cup K_2))= \Z * \Z$. This is a contradiction.

We summarize the above discussion with the following corollary. 

\cor \label{corollary:clifford-non-embedding}
	There exist pairs of totally real tori $K_1, K_2 \sub \R^4$ with the property that each $K_i$ is isotopic through totally real tori to the standard Clifford torus $\{(z_1, z_2) \in \C^2 \mid |z_1|=|z_2|=1 \}$, but such that $K_1 \cup K_2$ is not smoothly isotopic to a Lagrangian.
\ecor

\Cref{corollary:clifford-non-embedding}, and the example preceding it, can be generalized to $\bb{CP}^2$ and $S^2 \tms S^2$ by embedding $K_1, K_2$ into a polydisk contained in $\bb{CP}^2$ or $S^2 \tms S^2$. It can also be generalized to links with arbitrarily many components. 

\subsubsection{Linking classes}

A \emph{line field} on $\bb{T}^2$ is a $1$-dimensional sub-bundle of the tangent bundle $T \bb{T}^2$. There is an obvious correspondence between co-orientable line fields and nonvanishing $1$-forms on $\bb{T}^2$, up to scaling. Following Eliashberg--Polterovich \cite{eli-pol}, we say that a line field is \emph{trivial} if it agrees up to homotopy with the kernel of a linear $1$-form (this does not depend on the identification $\bb{T}^2= \R^2/ \bb{Z}^2$). 

Let $(M, \o)$ be a symplectic $4$-manifold and let $L \sub M$ be a nullhomologous Lagrangian torus. The map $TM|_L \to T^*L$ sending $v \mapsto \o(v,-)$ descends to an isomorphism 
\eq \label{equation:normal-bundle} TM|_L/ TL = NL \to T^*L,\eeq
where $NL \to L$ is the normal bundle. We fix an embedding $j:NL \to M$ with $j(0_L)=L$ furnished by the tubular neighborhood theorem.

A \emph{framing} of the normal bundle $NL \to L$ is a non-vanishing section. We say that two framings $\tau, \tau'$ are equivalent and write $\tau \sim \tau'$ if they are homotopic. Given a pair of framings $\tau_1, \tau_2$, there is a section $\phi:L \to \op{Aut}(NL)$ such that $\tau_1=\phi \tau_2$.  Passing to the quotient under $\sim$, we get a well-defined element $\delta(\tau_1, \tau_2)= [\phi] \in [L, \op{GL}_2^+(\R)]= [L, S^1]= H^1(L; \Z)$ called the \emph{difference class}. 

The nullhomologous Lagrangian torus $L \sub M$ inherits two canonical framings (up to homotopy). The \emph{nullhomologous framing} $\tau_N: L \to NL$ is determined by the property that the pushoff of $0_L$ in the direction of $\tau_N$ is mapped by $j$ to a nullhomologous torus in $M-L$. The \emph{Lagrangian framing} $\tau_L: L \to NL$ corresponds to the graph of a linear $1$-form on $\bb{T}^2$ under the isomorphism \eqref{equation:normal-bundle}.

\defi[see Sec.\ 1.1 in \cite{eli-pol}] Let $\s(L): H_1(L; \Z) \to \Z$ be defined by setting $\s(L)([a])= \op{lk}(L, a')$, where $[a] \in H_1(L; \Z)$, $a'$ is the pushoff of $a$ in the direction specified by the Lagrangian framing and $\op{lk}(-,-)$ is the linking number. We view $\s(L)$ as an element of $H^1(L; \Z)$ via the canonical isomorphism $H^1(L; \Z) = \op{Hom}_{\Z}(H_1(L;\Z), \Z)$ and call $\s(L)$ the \emph{(Polterovich) linking class} of $L \sub (M, \o)$. It is straightforward to show (see \cite[Sec.\ 3]{fintushel-stern}) that $\s(L)= \delta(\tau_N, \tau_L)$. 
\edefi

Eliashberg and Polterovich proved that $\s(L)=0$ for any Lagrangian torus $L \sub (\R^4, \o)$; see \cite[Thm.\ 1.1.A]{eli-pol}. The following corollary of \Cref{theorem:lagrangian-rational-introduction}, which is proved in \Cref{subsection:unlinking-rational} generalizes this to rational symplectic $4$-manifolds. 

\cor \label{corollary:vanishing-linking-class-rational} Let $(M, \o)$ be a rational symplectic $4$-manifold and let $L \sub M$ be a Lagrangian torus (which is nullhomologous by \Cref{proposition:all-tori-isotopic}). Then $\s(L)=0$. \ecor
 
We remark that \Cref{corollary:vanishing-linking-class-rational} cannot be generalized to arbitrary symplectic $4$-manifolds. Indeed, Fintushel and Stern \cite{fintushel-stern} constructed examples of closed symplectic $4$-manifolds containing nullhomologous Lagrangian tori with nonvanishing linking class. In fact, they also showed that these tori are in general not pairwise isotopic, in contrast to \Cref{proposition:all-tori-isotopic}.

As an application of \Cref{corollary:vanishing-linking-class-rational}, we have the following result, which was proved by Eliashberg and Polterovich for $M=\R^4$. 

\cor Let $(M, \o)$ be a rational symplectic $4$-manifold. Let $\S \sub M$ be an embedded closed $3$-manifold, whose characteristic foliation admits an embedded, invariant $2$-torus $L \sub \S$. If $L$ divides $\S$, then the characteristic line field restricted to $L$ is trivial. \ecor
\pf
Choose a vector field $v \in TM|_L$ which is tangent to $\S$ and transverse to $TL$ (i.e. a normal vector field to the embedding $L \sub \S$). Observe that $v$ defines a non-vanishing section of $NL= TM|_L/TL$ which is easily seen to be the nullhomotopic framing. Observe also that the characteristic line field on $L$ is precisely $\op{ker} \o(v, -)|_L$. 

Since $\s(L)=0$ by \Cref{corollary:vanishing-linking-class-rational}, the Lagrangian framing is equal to the nullhomologous framing. Hence $v$ corresponds to a linear $1$-form under the isomorphism \eqref{equation:normal-bundle}. But this is precisely what it means for $\op{ker} \o(v, -)$ to be a trivial line field.
\epf

\subsection*{Acknowledgements} We thank Yasha Eliashberg for helpful conversations. Part of this work was carried out when the second author visited the Department of Mathematics at Stanford University in February 2019, and when the first author visited the Department of Mathematics at Uppsala University in December 2019.

\section{Background material} \label{section:background-material}

This paper relies heavily on the analysis of pseudoholomorphic curves in symplectic manifolds and we assume that the reader has basic familiarity with the theory. As is customary in the literature, we will often drop the prefix ``pseudo" when discussing pseudoholomorphic curves, if it is clear from the context that we are not imposing any integrability assumptions. Similarly, we will sometimes blur the distinction between ``almost-complex structures" and ``complex structures". 

Given a Riemannian manifold $(L, g)$ let 
\begin{align*}
S^*_{r, g}(L) = \{ v \in T^* L \mid \|v\|_g = r \}, &\hspace{0.5cm} S^*_{<r,g}(L)= \{ v \in T^*L \mid \|v \|_g <r \}.
\end{align*}

We often write $S^*_rL$ and $S^*_{<r}L$ when the choice of metric is clear from the context. We also define $S^*_{\leq r, g} L$ in the obvious way.

We will make repeated use in this paper of \emph{neck-stretching} arguments, which were first introduced in symplectic topology in \cite{sftcompactness}. The basic principle of neck-stretching is to study limits of pseudoholomorphic curves in a symplectic manifold which degenerate near a Lagrangian submanifold. Neck-stretching arguments have become quite widespread in symplectic topology. We recall here the basic constructions for the purpose of fixing our notation. We refer the reader to \cite[Sec.\ 1.3]{sft} and \cite[Sec.\ 5]{evans} for a more detailed introduction. 

We begin with the following definition.
\defi \label{definition:neck-stretch-datum}
Fix a Lagrangian submanifold $L$ in a symplectic manifold $(M, \o)$. A \emph{neck-stretching datum} $ \eta \in \mc{D}(M, \o; L)$ consists of the following objects (see \cite[Sec.\ 2.5]{dgi}): 
\begin{itemize}
\item A metric $g$ on $L$.
\item A Weinstein embedding $\phi: (N(0_{L}), 0_{\bb{T}^2}) \hookrightarrow (M, L)$, where $N(0_{\bb{T}^2}) \sub T^* L$ is an open neighborhood of $0_{L}$ which contains $S^*_{\leq 4, g}L$.
\item An almost-complex structure $J_{\op{cyl}}$ on $\R \tms S^*_{1,g}L$ which is compatible with $d\l_{can}$ and translation-invariant.
\item An almost-complex structure $J_{\op{std}}$ on $T^* L$ which is compatible with $d \l_{can}$ and coincides with $J_{\op{cyl}}$ outside $T^*_2L$.
\item An almost-complex structure $J_{\infty}$ on $M-L$ which is compatible with $\o$ and whose pullback under $\phi$ coincides with $J_{\op{cyl}}$ on $S^*_{\leq 4, g}L- 0_{L}$. 
\end{itemize}
\edefi

For any given $L \sub (M, \o)$, there is always an abundance of choices of neck-stretching data. However, it's often convenient to make a particular choices of datum for the purposes of analyzing the SFT limit discussed below.  In this paper, we will only ever consider neck-stretching data associated to Lagrangian tori. In this case, one construction of the relevant almost-complex structures is described in \cite[Sec.\ 4]{dgi}. 

Given a choice of neck-stretching datum, the \emph{neck-stretching procedure} produces a family $J_{\tau}$, $\tau \geq 0$, of $\o$-compatible almost-complex structures on $(M, \o)$. This procedure is described in \cite[Sec.\ 2.5]{dgi}, \cite[Sec.\ 3.4]{sftcompactness}. Given a sequence of $J_{\tau_i}$-holomorphic curves $u_i$, for $\tau_i \to \infty$, the SFT compactness theorem ensures that (after possibly passing to a subsequence) the $u_i$ converge to a \emph{holomorphic building} $\mathbf{u}$. 

\defi 
A holomorphic building $\mathbf{u}$ associated to a neck-streching datum $\eta \in \mc{D}(M, \o; L)$ consists in the following data:
\begin{itemize}
\item A nodal curve $C$.
\item A collection $\mc{C}_1$ of punctured $J_{\infty}$-holomorphic curves in $M-L$. These are said to be \emph{top level} curves.
\item A collection $\mc{C}_2$ of punctured $J_{\op{cyl}}$-holomorphic curves in $\R \tms S^*_1L$. These are said to be \emph{middle level} curves.
\item A collection $\mc{C}_3$ of punctured $J_{\op{std}}$-holomorphic curves in $T^*L$. These are said to be \emph{bottom level} curves.
\end{itemize}

The elements of $\mc{C}_1 \cup \mc{C}_2 \cup \mc{C}_3$ are required to be in bijective correspondence with the components of $C$. Moreover, any pair of components of $C$ which share a node must correspond to a pair of punctured curves on adjacent levels which share a puncture.  If $C$ is a nodal sphere, we say that $\mathbf{u}$ is a \emph{split} (or \emph{broken}) sphere. We refer to \cite[Sec.\ 2.4]{dgi} or \cite[Sec.\ 9]{sftcompactness} for a more detailed definition. 
\edefi

A holomorphic building $\bf{u}$ can be compactified to form a continuously embedded surface $\S$ inside $M$.  We then say that $\bf{u}$ represents the class $[\S] \in H_2(M; \Z)$. If $\bf{u}$ is a limit of holomorphic curves of class $\a \in H_2(M; \Z)$, then it can be shown that $[\S] = \a$.  

When we write informally that we \emph{stretch the neck} along a Lagrangian $L \sub (M, \o)$, this means that we apply the neck-stretching procedure for some choice (possibly unspecified) of neck-stretching datum.  If $L$ is a finite disjoint union of connected Lagrangian submanifolds $L_1,\dots,L_n$, then we may apply the neck-stretching procedure to the disconnected Lagrangian $L= \cup_{i=1}^n L_i$. In this case, we say informally that we stretch the neck simultaneously along the $L_i$.

A pseudoholomorphic curve with domain $\bb{CP}^1 - \{\infty\}= \C$ is called a pseudoholomorphic \emph{plane}, while a pseudoholomorphic curve with domain $\bb{CP}^1 - \{0, \infty\}$ will be called a pseudoholomorphic \emph{cylinder}.

In this paper, we will typically be considering neck-stretching data associated to a Lagrangians  $L \sub (S^2 \tms S^2, \o \oplus \o)$. Following \cite{dgi}, a split sphere of class $[S^2 \tms *]$ or $[* \tms S^2]$ is said to be of \emph{Type I} if it consists of a single bottom level cylinder along with two top-level planes, each of which is asymptotic to a geodesic of Maslov index $2$. Otherwise, the split sphere is said to be of \emph{Type II}. In case $L$ is disconnected, observe that the components of a Type I building have punctures asymptotic \emph{a single} component of $L$.

\section{Hamiltonian unlinkedness of weakly exact tori} \label{section:weakly-exact}

The goal of this section is to prove \Cref{thm:exact}, which was stated as \Cref{theorem:weakly-exact-classification} in the introduction. This theorem gives a classification up to Hamiltonian isotopy of weakly exact, rational Lagrangian tori in $T^*\TT^2 - 0_{\TT^2}$. We also prove \Cref{corollary:eli-pol-giroux} and \Cref{corollary:pi-1-cotangent}, which are both straightforward consequences of the theorem. 

\subsection{Some background} \label{subsection:some-background}
Let us fix an identification \eq \label{equation:torus-identification} \bb{T}^2 = \R^2 / \Z^2.\eeq This induces an indentification $T^* \bb{T}^2 = \bb{T}^2 \tms \R^2$ with coordinates $(q_1, q_2, p_1, p_2) \in \R^2 / \Z^2 \tms \R^2$. We let $\l_{\op{can}} = p\,dq= p_1 dq_1+ p_2 dq_2$ be the tautological $1$-form on $T^* \bb{T}^2$.


Recall the following definitions. 

\defi \label{definition:weaklyexact}  Given a symplectic manifold $(M,\o)$ and a Lagrangian $L \sub M$, we say that $L$ is \emph{weakly exact} if $\omega(u)= \int_{D^2} u^* \omega=0$ for all maps $u: (D^2, \d D^2) \to (M, L)$. \edefi 

\defi \label{definition:rational}  A Lagrangian torus $L \subset (T^*\TT^2,d(p\,dq))$ is said to be \emph{rational} if the symplectic action class $p\,dq|_{L}$ is contained in the subset $c\cdot H^1(L,\Q) \sub H^1(L; \R)$ for some real number $c>0$. \edefi 

In this section we consider a weakly exact Lagrangian torus $L \subset (T^*\TT^2,p\,dq)$ which is disjoint from the zero section $0_{\TT^2} \subset T^*\TT^2.$

\begin{theorem} \label{thm:exact}
Let $L\subset T^*\TT^2$ be a weakly-exact, rational Lagrangian torus. If $L$ is disjoint from the zero section, then $L$ is Hamiltonian isotopic \emph{in the complement of the zero section} to a standard torus $\TT^2 \times \{\mathbf{p}^L\}$, where $\mathbf{p}^L \in \R^2- \{0\}$. 
\end{theorem}

Since $L$ is weakly exact, it's not hard to show that $L$ becomes exact after translation by the graph of a suitable closed $1$-form; see \cite[Lem.\ 9.2]{dim1}. It then follows from the nearby Lagrangian conjecture for $\bb{T}^2$ proved in \cite[Thm.\ B]{dgi} that $L$ is Hamiltonian isotopic to a standard torus $\TT^2 \times \{\mathbf{p}^L\}$. The content of \Cref{thm:exact} is thus that the isotopy can be confined to the complement of the zero section.

\rmk It follows from the classical non-displaceability results of Laudenbach--Sikorav \cite{lau-sik} or Floer \cite{floer} that a torus which satisfies the hypotheses of \Cref{thm:exact} cannot be exact. 
\ermk

\rmk One could strengthen \Cref{thm:exact} to include Lagrangian tori which are weakly exact inside the smaller symplectic manifold $T^*\TT^2 - 0_{\TT^2}$. In fact, it can be shown by a neck-stretching analysis that such tori are automatically weakly exact in $T^*\TT^2$.
\ermk

The heart of the proof of \Cref{thm:exact} consists in establishing the following proposition.

\begin{proposition}
\label{prp:IntoConvex}
Suppose that $L\subset T^*\TT^2$ is a weakly exact Lagrangian torus which satisfies $L \cap 0_{\TT^2}=\emptyset$. Then there exists a convex subset $U \sub \R^2 - \{0 \}$ and a Lagrangian isotopy
\begin{gather*}
\varphi_t \colon \TT^2 \to L^t \subset T^*\TT^2 - 0_{\TT^2},\\
L^0=L, \:\:\: L^1 \subset \TT^2 \times U,
\end{gather*}
for $t \in [0,1]$. If $L$ is rational, then we may assume that the symplectic action satisfies
$$[\varphi_t^*(p\,dq)]=e^{g(t)}\cdot[\varphi_0^*(p\,dq)] \in H^1(\TT^2;\R)$$
for some smooth path $g: [0,1] \to \R$. 
\end{proposition}

We will now show that \Cref{thm:exact} is an easy consequence of the above result combined with \cite[Thm.\ B]{dim1}. We begin by recalling the following facts, the first of which is standard. 

\lem \label{lemma:lag-ham}
	Let $(M, d \l)$ be an exact symplectic manifold, and let $\phi_s: L \to (M, d \l)$ be a Lagrangian isotopy, for $s \in [0,1]$. This isotopy is Hamiltonian if and only if $[\phi_s^*(\l)] \in H^1(L; \R)$ is constant (i.e. independent of $s$).
\elem
\pf
	Let $\a_s:= \d_s \phi^*_s\l \in \O^1(L)$. Letting $\Phi: [0,1] \tms L  \to M$ be defined as $\Phi(s, x):= \phi_s(x)$, one computes that $\Phi^*(\o)= \Phi^*(d\l)= d\Phi^*(\l)= ds \wedge \a_s$.  It now follows from \cite[Exercise 6.1]{polt-book} that the isotopy $\phi_s$ is Hamiltonian if and only if $\a_s$ is exact for all $s \in [0,1]$. But $\a_s$ is exact if and only if $[\phi_s^*(\l)]$ is constant. \epf

\lem[see Lem.\ 9.1 in \cite{dim1}]  \label{lemma:weakly-exact}
Let $L \sub T^* \bb{T}^2$ be a Lagrangian torus. The following are equivalent: 
\begin{itemize}
\item $L$ is weakly exact,
\item $L$ is homologically essential in $ T^* \bb{T}^2$,
\item the inclusion $L \hookrightarrow T^* \bb{T}^2$ is a homotopy equivalence.
\end{itemize}
\elem

\thm[Dimitroglou Rizell, see Thm.\ B in \cite{dim1}] \label{theorem:1-form-dim} Let $U \sub \R^2$ be a convex subset and let $L \sub \bb{T}^2 \tms U \sub T^*\bb{T}^2$ be a weakly exact Lagrangian torus. Then $L$ is Hamiltonian isotopic to the graph of a closed $1$-form via an isotopy which is supported in $\bb{T}^2 \tms U$. \ethm

Our identification $\bb{T}^2= \R^2 / \Z^2$ in \eqref{equation:torus-identification} induces a basis 
\eq \label{equation:splitting}\langle \mathbf{e}_1,\mathbf{e}_2 \rangle = H_1(\bb{T}^2; \R). \eeq 
If $L \sub T^*\bb{T}^2$ is weakly exact, \Cref{lemma:weakly-exact} implies that the composition  $L \hookrightarrow T^*\bb{T}^2 \to \bb{T}^2$ is a homotopy equivalence. We therefore obtain a basis $$\langle \mathbf{e}^L_1,\mathbf{e}^L_2 \rangle=H_1(L)$$ by pulling back $\mathbf{e}_1$ and $\mathbf{e}_2$. We set $$p^L_i \coloneqq \int_{\mathbf{e}^L_i}p\,dq \in \R$$ and write $\mathbf{p}^L = (p^L_1,p^L_2) \in \R^2$.

\begin{proof}[Proof of \Cref{thm:exact}] 

The path of tori $\tilde{L}^t\coloneqq e^{-g(t)}\cdot L^t$ given as the image of $L^t$ under a path of suitable fiberwise rescalings is generated by a Hamiltonian isotopy by \Cref{lemma:lag-ham}. This path thus defines a Hamiltonian isotopy which takes $\tilde{L}^0=L$ into $\tilde{L}^1 \subset \TT^2 \times (e^{-g(1)}\cdot U)$ where $e^{-g(1) }\cdot U\subset \R^2 - \{0\}$ is again convex. Let us set $U_1:= e^{-g(1) }\cdot U$. 

According to \Cref{theorem:1-form-dim}, $\tilde{L}^1$ isotopic to the graph of a closed $1$-form $\a^1 \in \O^1(\bb{T}^2)$, via an isotopy which is supported in $\bb{T}^2 \tms U_1$.  We write $\a= f(q_1,q_2) dq_1 + g(q_1, q_2) dq_2$ and observe that $(f(q_1, q_2), g(q_1, q_2)) \in U_1$ for all $(q_1, q_2) \in \bb{T}^2$. 

We have that
\begin{align*} (p_1^{\tilde{L}^1}, p_2^{\tilde{L}^2}) &= (\l(\mathbf{e}_1^{\tilde{L}^1}), \l(\mathbf{e}_2^{\tilde{L}^1}) )  \\
&=(p_1 dq_1(\mathbf{e}_1^{\tilde{L}^1}), p_2 dq_2(\mathbf{e}_2^{\tilde{L}^1}) ) \nonumber \\
&=(p_1 dq_1( \mathbf{e}_1^{\tilde{L}^1}+ \mathbf{e}_2^{\tilde{L}^1}), p_2 dq_2( \mathbf{e}_1^{\tilde{L}^1}+ \mathbf{e}_2^{\tilde{L}^1}) ) \nonumber \\
&= \int_0^1 (f(t, t), g(t,t)) dt \in U_1,
\end{align*} where the last line uses the fact that $U_1$ is convex. On the other hand, since $\tilde{L}^t$ is a Hamiltonian isotopy, it follows from \Cref{lemma:lag-ham} that we have $(p_1^{\tilde{L}^1}, p_2^{\tilde{L}^2}) = (p_1^L, p_2^L)= \mathbf{p}^L$. 

Let $\a^0= p_1^L dq_1+ p_2^L dq_2$ and let $\a^s = (1-s) \a^0 + t \a^s$ for $s \in [0,1]$. Then the family $\hat{L}^s= \op{graph}(\a^s)$ is contained in $\bb{T}^2 \tms U_1$ since $U_1$ is convex. It follows by \Cref{lemma:lag-ham} that $\hat{L}^s$ is a Hamiltonian isotopy, which completes the proof.
\end{proof}

\pf[Proof of \Cref{corollary:eli-pol-giroux}]
Given any Weinstein neighborhood $\mc{U}$ of $L$ (which can be assumed to be disjoint from the zero section), $L$ is Lagrangian isotopic to a rational torus $L' \sub \mc{U}$. It then follows from \Cref{lemma:weakly-exact} and \Cref{thm:exact} that $L'$ is Lagrangian isotopic to $\bb{T}^2 \tms \{\mathbf{p} \}$, for some $\mathbf{p} \neq (0,0)$. 
\epf

\pf[Proof of \Cref{corollary:pi-1-cotangent}]
Given a closed $1$-form $\a \in \G(T^*\bb{T}^2)$, let $\tau_{\a}: T^*\bb{T}^2 \to T^*\bb{T}^2$ be the symplectomorphism defined by translation by the graph of $\a$. According to \cite[Lem.\ 9.2]{dim1}, there exists a suitable $1$-form $\a_1 \in \G(T^*\bb{T}^2)$ such that $\tau_{\a_1}L_1$ is exact. Hence, after replacing the $L_i$ by $\tau_{\a_1}(L_i)$ for $i=1,2$, we may assume that $L_1$ is exact. By the Nearby Lagrangian Conjecture for $\bb{T}^2$ (see \cite[Thm.\ 7.1]{dgi}), $L_1$ is Hamiltonian isotopic to the zero section. So we may as well assume that $L_1=0_{\bb{T}^2}$. The claim now follows by \Cref{corollary:eli-pol-giroux}. 
\epf

The remainder of \Cref{section:weakly-exact} is devoted to proving \Cref{prp:IntoConvex}.

\subsection{Some notation and reductions} \label{subsection:notation-reductions} 

We begin with a sequence of reductions.

Observe first that any Lagrangian torus $L \sub T^*\bb{T}^2- 0_{\bb{T}^2}$ admits a Lagrangian isotopy to a rational Lagrangian torus $L'$ which is contained in an arbitrarily small Weinstein neighborhood of $L$. It follows from \Cref{lemma:weakly-exact} that $L'$ is weakly-exact if and only if $L$ is weakly exact. For the purpose of proving \Cref{prp:IntoConvex}, there is therefore no loss of generality in assuming that $L$ is rational. 

Next, it is clear from the statement of \Cref{prp:IntoConvex} that we can freely rescale the fibers of $T^* \bb{T}^2$ by a constant $\a>0$. In particular, we may assume that $L$ is \emph{integral}, i.e. $p dq|_L \in c \dot H^1(L; \Z)$ for some real number $c>0$. Similarly, we can freely act by $\op{SL}(2,\Z)$ on $T^* \bb{T}^2$, where the action is induced by the natural action on $\bb{T}^2 = \R^2/ \Z^2$. 

By applying a suitable element of $\op{SL}(2,\Z)$, we can therefore assume in the statement of \Cref{prp:IntoConvex} that $(p^L_1,p^L_2) =(0, c)$. After possibly applying a further fiberwise recaling, we can further assume that $c \in (0,1/4)$ and that $L$ is contained in the disk bundle $D^*_{1/4} \bb{T}^2$ consisting of covectors of norm less than $1/4$ with respect to the metric on $\bb{T}^2$ induced by \eqref{equation:torus-identification}. These assumptions will be in effect for the remainder of \Cref{section:weakly-exact}.

Let us endow the sphere $S^2$ with the standard triple $(\o, i, g_{S^2})$, where $\o$ is the Fubini-Study form, $i$ is the standard integrable complex structure and $g_{S^2}$ is the round metric. We normalize these so that $\o= \op{dvol}_{g_{S^2}}$ has volume $1$. Let $S^1 \sub S^2$ be the equator and let $\{0\}$ and $\{\infty \}$ be the south pole and north pole, respectively. 

We fix an identification $S^1= \R/ \Z$ and endow $S^1$ with the standard Euclidean metric. We let $D^*_r S^1$ be the disk bundle of covectors of norm $<r$, measured with respect to this metric. 

Fix an exact symplectomorphism
\eq \label{equation:weinstein-s1} D^*_{1/2}S^1 \xhookrightarrow{\simeq} S^2 - \{0\} - \{\infty\}, \eeq 
which takes $S^1$ to the equator and takes the fibers to geodesic arcs connecting $\{0\}$ and $\{\infty \}$.  

We now consider the product manifold $S^2 \tms S^2$ endowed with the symplectic form $\o \oplus \o$. Let $$D_\infty = S^2 \times \{0,\infty\} \:\: \cup \:\: \{0,\infty\} \times S^2$$ be the nodal divisor consisting of four lines. Observe that \eqref{equation:weinstein-s1} induces a Weinstein embedding 
\eq \label{equation:weinstein-t2} D^*_{\a} \bb{T}^2 \sub D^*_{1/2}S^1 \tms D^*_{1/2}S^1 \xhookrightarrow{\simeq} S^2 \tms S^2 - D_{\infty}, \eeq 
where the second map is an exact symplectomorphism and $\a< 1/2$. 

We let $L_0$ denote the image of the zero section along the embedding \eqref{equation:weinstein-t2}, which is identified with the product of equators $S^1 \tms S^1 \sub S^2 \tms S^2$. Similarly, we let $$L_1 \sub S^2 \times S^2 - ((S^1 \times S^1) \cup D_\infty)$$ denote the image of $L$ under the embedding \eqref{equation:weinstein-t2}. This notation will be in effect for the remainder of this section.

\subsection{Preliminary neck-stretching analysis}  \label{subsection:preliminary-neck-stretch}
Let $J_{\op{cyl}}$ be the standard cylindrical almost-complex structure on $T^*\TT^2 - 0_{\TT^2}$. In our coordinates $(q_1, q_2, p_1, p_2)$, it is defined by the equation $J_{\op{cyl}} \d_{q_i} = - \|(p_1, p_2) \| \d_{p_i}$, where $\|- \|$ is the Euclidean metric (cf.\ \cite[Sec.\ 4]{dgi}). 

For the remainder of this section, we fix a neck-stretching datum $\eta= (g, \phi, J_{\op{cyl}}, J_{\op{std}}, J_{\infty}) \in \mc{D}(S^2 \tms S^2, \o \oplus \o; L_0 \cup L_1)$ satisfying the following properties. 
\begin{enumerate}[label=(A\arabic*)]
\item \label{item:metric-flat} The metric $g$ on $L_0 \cup L_1$ is flat. 
\item \label{item:disjoint-image} The image of $\phi$ is disjoint from $D_{\infty}$. Letting $\phi_0$ (resp.\ $\phi_1$) denote the restriction of $\phi$ to the cotangent bundle of $L_0$ (resp.\ $L_1$), we have that $\phi_0$ coincides with the embedding \eqref{equation:weinstein-t2}, for $\a>0$ suitably small. 
\item \label{item:standard-ac} $J_0= i$ in some fixed neighbourhood $\mc{U}_{\infty}$ of $D_\infty$, such that $\mc{U}_{\infty}$ is disjoint of $\op{Im} \phi$. 
\item \label{item:standard-lines} For $\t \in S^1$, the standard lines $S^2 \times \{e^{i\theta}\} - L_0$ are $J_{\infty}$-holomorphic near $L_0$. 
\item $J_{\op{cyl}}$ and $J_{\op{std}}$ are almost-complex structures on $R \tms S^*_{1,g}\bb{T}^2$ and $T^* \bb{T}^2$ respectively, defined as in \cite[Sec.\ 4]{dgi}.\footnote{Strictly speaking, to be consistent with \Cref{definition:neck-stretch-datum}, we should view $J_{\op{cyl}}$ and $J_{\op{std}}$ as being defined in the cotangent bundle of the disjoint union $\bb{T}^2 \sqcup \bb{T}^2$. We hope the reader will forgive this slight abuse of notation.}
\end{enumerate}

It is straightforward to check that such a choice of neck-stretching data exists. 


It will be convenient to record the following lemma, which is an immediate consequence of \ref{item:disjoint-image} and the definition of $J_{\op{std}}$ (see also \cite[Lem.\ 4.2]{dgi} for a description of all $J_{\op{std}}$-holomorphic cylinders of finite energy). 
\lem \label{lemma:standard-cylinders}
	The cylinders $u_{\t}:= \{q_2= \t, p_2=0\} \sub T^*S^1 \tms T^*S^1= T^* \bb{T}^2$ are $J_{\op{cyl}}$-holomorphic. The restriction of $u_{\t}$ to the domain of $\phi_0$ is mapped into the line $S^2 \tms \{e^{i \t} \}$.
\elem
\qed

Let $A=  [S^2 \times \{\pt\}] \in H_2(S^2 \times S^2)$.  Since $A$ is a primitive class of minimal area, Gromov's celebrated result (see \cite[0.2.A]{gromov}) ensures that there is a foliation of $S^2 \tms S^2$ by $J_{\tau}$-holomorphic lines in the class $A$, for all $\tau \geq 0$. The following lemma describes the split spheres which may result from applying the SFT compactness theorem to $J_{\tau_i}$-holomorphic spheres in class $A$, for some sequence $\tau_i \to \infty$. 

\begin{figure}[htp]
\centering
\vspace{3mm}
\hspace{20mm}
\labellist
\pinlabel $S^2\times S^2-(L_0\cup L_1)$ at -48 58
\pinlabel $T^*L_i$ at -18 19
\pinlabel $D_\infty$ at 70 81
\pinlabel $D_\infty$ at 29 81
\pinlabel $\color{red}\mathbf{e}^{L_i}_1$ at 30 22
\pinlabel $\color{red}-\mathbf{e}^{L_i}_1$ at 67 22
\pinlabel $D_\infty$ at 143 81
\pinlabel $D_\infty$ at 264 81
\pinlabel $\color{red}\mathbf{e}^{L_i}_1$ at 145 21
\pinlabel $\color{red}-\mathbf{e}^{L_i}_1$ at 181 21
\pinlabel $\color{red}-\mathbf{e}^{L_j}_1$ at 224 21
\pinlabel $\color{red}\mathbf{e}^{L_j}_1$ at 266 21
\endlabellist
\includegraphics{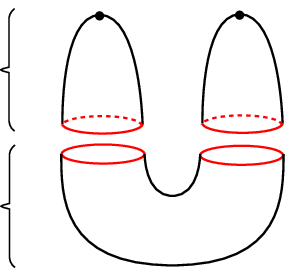}
\hspace{15mm}
\includegraphics{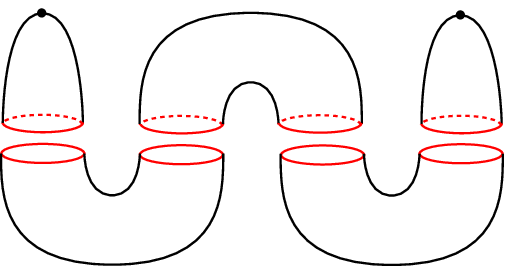}
\caption{On the left: a split sphere of Type I (the generic configuration), with each of the planes intersecting $D_\infty$ transversely in a single point. On the right: a split sphere of Type II (the exceptional configuration). The Type II configuration is shown in \Cref{lma:notypeII} not to occur.} 
\label{fig:breaking}
\end{figure}

\begin{lemma}
  \label{lma:notypeII}
Any split sphere representing the class $A$ consists in at most two components in its top level. Moreover, each such component must be a plane asymptotic to a geodesic in one of the classes $\pm \mathbf{e}^{L_i}_1$. In particular, there are no Type II configurations (see \Cref{fig:breaking}). 
\end{lemma}

\begin{proof}
Positivity of intersection, together with the holomorphicity of $D_\infty$, shows that any pseudoholomorphic line (broken or not) in the homology class $A \in H_2(S^2 \times S^2)$ is disjoint from $S^2 \times \{0,\infty\}$, and intersects each component of $\{0,\infty\} \times S^2$ transversely in a unique point.

Since $L_0$ and $L_1$ were endowed with the flat metric, it follows that any closed geodesic on either torus is homologically essential inside $S^2 \times S^2 - D_\infty$. (Here we have used the assumption that $L_0$ and $L_1$ both are homologically essential.) Hence all pseudoholomorphic planes asymptotic to either torus must intersect $D_\infty$. Together with the intersection properties established in the previous paragraph, we conclude that there exists precisely two planes and, moreover, that any plane asymptotic to $L_1$ (resp. $L_0$) is asymptotic to geodesics in the classes $\pm \mathbf{e}^{L_1}_1$ (resp. $\pm \mathbf{e}_1^{L_0}$). 

The remaining top level components of the building must thus consist of cylinders contained inside $S^2 \times S^2 - (L_0 \cup L_1 \cup D_\infty)$. An elementary topological argument, which takes the asymptotics of the planes as established above into account, now shows that any puncture of a sphere that arises in a broken line in homology class $A$ can be asymptotic to geodesics only in homology classes of the form $\pm \mathbf{e}_1^{L_0}$ or $\pm \mathbf{e}_1^{L_1}$. (Here we have again used the assumption that $L_0$ and $L_1$ both are homologically essential.)

The statement finally follows from the topological fact that any cylinder as above has vanishing symplectic area, by the assumptions on the symplectic action classes of $L_i$; hence it cannot be pseudoholomorphic for a compatible almost-complex structure.
\end{proof}

\subsection{Constructing a fibration}
We now collect some auxiliary results for constructing and manipulating symplectic $S^2$-fibrations. We will put together these results in \Cref{subsection:finalproof} to prove \Cref{prp:IntoConvex}. The reader may wish to pass directly to \Cref{subsection:finalproof} and consult this section as the need arises.

\begin{figure}[htp]
\vspace{3mm}
	\labellist
        \pinlabel $y$ at 153 47
        \pinlabel $0$ at 153 80
        \pinlabel $\infty$ at 163 139
        \pinlabel $x$ at 87 53
        \pinlabel $\color{blue}L_0$ at 60 108
        \pinlabel $\color{blue}L_1$ at 195 106
        \pinlabel $\pi_1^{-1}(y)$ at 195 128
        \pinlabel $\pi_1^{-1}(x)$ at 55 124
	\pinlabel $0$ at -5 33
        \pinlabel $\infty$ at 94 142
        \pinlabel $0$ at 94 96
	\pinlabel $\infty$ at 213 33
	\pinlabel $\color{blue}\pi_1(L_1)$ at 135 33
	\pinlabel $\color{blue}\pi_1(L_0)$ at 75 33
	\endlabellist
	\includegraphics[scale=1]{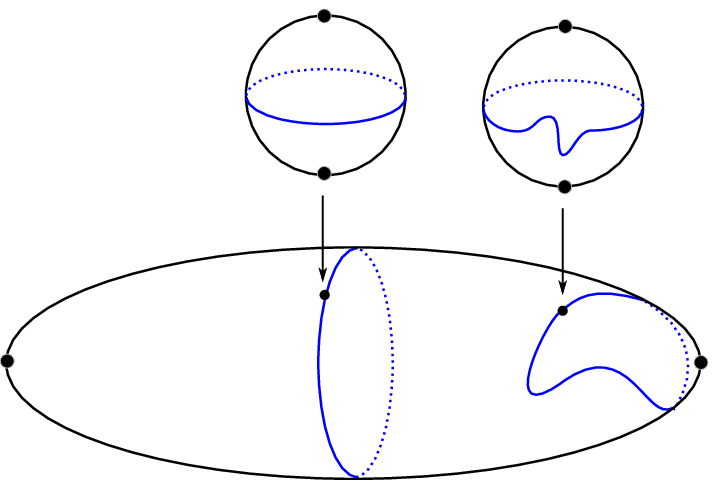}
	\caption{A fibration $\pi_1$ compatible with both tori $L_i$, $i=0,1$.}
	\label{fig:compatible}
\end{figure}

\begin{proposition}
  \label{prp:fibrations}
Let $L_1 \subset S^2 \times S^2 - D_\infty$ be a Lagrangian torus that is disjoint from $L_0$. Then there exists a smooth family $J^t$ of compatible almost-complex structures on $(S^2 \tms S^2, \o \oplus \o)$, and a family of smooth symplectic $S^2$-fibrations
  $$\pi_t \colon S^2 \times S^2 \to S^2, \:\: t \in [0,1]$$
with $J^t$-holomorphic fibers.
The complex structures satisfy the following properties:
\begin{itemize}
\item[(i)] $J^t$ is independent of $t$ near $\{0, 1\}$, 
\item[(ii)] $J^t=i$ near $D_{\infty}$. 
\end{itemize}
The fibrations $\pi_t$ satisfy the following properties: 
\begin{enumerate}
\item \label{item:t=0} $\pi_0^{-1}(S^1)=\op{pr}_2^{-1}(S^1) \supset L_0$, where $\op{pr}_2 \colon S^2 \times S^2 \to S^2$ is the canonical projection to the second factor.
\item \label{item:over-poles} $\pi_t^{-1}(\infty)=S^2 \times \{\infty\}$ and $\pi_t^{-1}(0)=S^2 \times \{0\}$, for all $t \in [0,1]$,
\item \label{item:over-equator} $\pi_t(L_0)=S^1 \subset S^2$ is the equator for all $t\in[0,1]$. Each fiber $\pi_t^{-1}(e^{i\theta})$ over a point $e^{i\theta} \in S^1$ in the equator coincides with a standard line $S^2 \times \{e^{if_t(\theta)}\}$ inside some small fixed neighbourhood of $L_0$, where $f_t \colon S^1 \to S^1$ is a family of diffeomorphisms that depend smoothly on $t \in [0,1]$.
\item \label{item:t=1} $\pi_1(L_1) \subset S^2 - \{\infty \cup 0 \cup S^1\}$ is a homotopically essential embedded closed curve that is contained in the complement of the equator; see \Cref{fig:compatible}.
\end{enumerate}
\end{proposition}

\pf  
Some arguments in this proof are drawn from the proof of \cite[Thm.\ D]{dgi} and will therefore only be summarized. We start by stretching the neck around $L_0$ and $L_1$ simultaneously using the neck-stretching datum $\eta=(g, \phi, J_{\op{cyl}}, J_{\op{std}}, J_{\infty})$ introduced in \Cref{subsection:preliminary-neck-stretch}. We then consider the SFT-limit of lines in class $A \in H_2(S^2 \times S^2)$. By \Cref{lma:notypeII}, we infer that any broken line is asymptotic to precisely one of the two tori $L_i$, for $i=0,1$. 

Let $\G^0(L_i)$ (resp.\ $\G^{\infty}(L_i)$) be the moduli space of $J_{\infty}$-holomorphic planes asymptotic to geodesics in the class $\pm \mathbf{e}^{L_i}_1$ which intersect $\{0\} \tms S^2$ (resp.\ $\{\infty\} \tms S^2$). The analysis of \cite[Sec.\ 5.2]{dgi} shows there is a unique plane in $\G^0(L_i)$ and in $\G^{\infty}(L_i)$ which is asymptotic to each geodesic $(L_i, g)$ in the class $\pm \mathbf{e}^{L_i}_1$ (recall from \ref{item:metric-flat} that $g$ is the flat metric). 

The smoothing procedure of \cite[Prop.\ 5.16]{dgi} allows us to glue together the two families of planes $\G^0(L_1)$ and $\G^{\infty}(L_1)$. We obtain a hypersurface $H^1$ which contains $L_1$, and an almost-complex structure $J^1_{\infty}$ on $S^2 \tms S^2-L_0$ such that such that $H^1$ is foliated by $J^1_{\infty}$-holomorphic spheres in the class $A$, and such that $J^1_{\infty}=J_{\infty}$ in $S^2 \tms S^2-(L_0 \cup \op{Im} \phi_1)$. 

Changing gears slightly, it follows from \ref{item:disjoint-image}, \ref{item:standard-ac} and \ref{item:standard-lines} that there exists a compatible, cylindrical almost-complex structure $J^0_{\infty}$ on $S^2 \tms S^2-L_0$ with the following properties:
\begin{itemize}
\item[(a)] $\phi_0^*J^0_{\infty}= J_{\op{cyl}}$ on $S^*_{\leq 4, g} \bb{T}^2- 0_{\bb{T}^2}$,
\item[(b)] $\op{pr}_2^{-1}(e^{i\t}) - L_0 $ is $J^0_{\infty}$-holomorphic, for $\t \in S^1$, 
\item[(c)] $J^0_{\infty}=i$ in $\mc{U}_{\infty}$. 
\end{itemize}

Let $J^t_{\infty}, t \in [0,1]$ be family of compatible almost-complex structures on $S^2 \tms S^2-L_0$ interpolating between $J^0_{\infty}$ and $J^1_{\infty}$. We may assume that $J^t_{\infty}$ is independent of $t$ near the endpoints, that $\phi_0^*J^0_{\infty}= J_{\op{cyl}}$ on $S^*_{\leq 4, g} \bb{T}^2- 0_{\bb{T}^2}$, and that $J^t_{\infty}=i$ in $\mc{U}_{\infty}$. 

Let $\G^0_t(L_0), \G^{\infty}_t(L_0)$ be the moduli space of $J^t_{\infty}$-holomorphic planes asymptotic to geodesics in the class $\pm \mathbf{e}^{L_0}_1$ which intersect $\{0\} \tms S^2$ (resp.\ $\{\infty\} \tms S^2$). These moduli spaces are compact by \Cref{lma:notypeII}, and Wendl's automatic transversality result \cite[Thm.\ 1]{wendlautomatic} implies that they all consist of regular curves. We now again apply the smoothing procedure of of \cite[Prop.\ 5.16]{dgi}. (Although the smoothing procedure of \cite[Prop.\ 5.16]{dgi} was only described for a single stretched almost-complex structure, it can also be performed with smooth dependence on a one parameter family of almost-complex structures without additional work.) We thus obtain a family of hypersurfaces $H^0_t$ which contain $L_0$, and a family of almost-complex structure $J^t$ on $S^2 \tms S^2$, such that $J^t= J^t_{\infty}$ in $S^2 \tms S^2-(\op{Im} \phi_0)$. The hypersurface $H^0_t$ is foliated by $J^t$-holomorphic spheres in the class $A$. Moreover, the proof of \cite[Prop.\ 5.16]{dgi} also ensures that these $J^t$-holomorphic spheres agree in a neighborhood of $L_0$ with the image under $\phi$ of the standard cylinders $u^{\t}$, for $\t \in S^1$, described in \Cref{lemma:standard-cylinders}. 

It follows from Gromov's celebrated result \cite[0.2.A]{gromov} that the foliation of $H^0_t$ extends to a $1$-parameter family of foliations of $S^2 \tms S^2$ by $J^t$-holomorphic spheres in the class $A$. We define a family of fibrations $\pi_t: S^2 \tms S^2$ in the usual way: given $p \in S^2 \tms S^2$, we let $\pi_t(p) \in S^2$ be the unique intersection point of a $J^t$-holomorphic sphere passing through $p$ with $\{0\} \tms S^2 \equiv S^2$.

It is clear by construction that (i) and (ii) are satisfied. \eqref{item:t=0} follows by the definition of $J^0$ and (a)-(c). \eqref{item:over-poles} follows from the fact that $J^t=i$ in $\mc{U}_{\infty}$. \eqref{item:over-equator} follows from the construction of $H^0_t$, as was already observed above. Finally, \eqref{item:t=1} follows from the fact that $H^0_1$ and $H^1$ are both foliated by $J^1$-holomorphic spheres in the class $A$, and must therefore be disjoint by positivity of intersection since $[A]^2=0$. \epf

In the following we let
$$\ell_1 \coloneqq S^2 \times \{1\},$$
be a standard holomorphic line in $S^2 \times S^2$ which intersects $L_0$ in a smooth circle. We identify $T_p \R^4$ with $\R^4$ for any $p \in \R^4$. If $\S \sub \R^4$ is a $d$-dimensional submanifold which is contained in a $d$-dimensional affine subspace of $\R^4$, we let $T\S$ be the tangent space at any point of $\S$. For $r>0$, we let $D_r \sub \R^2$ be the closed disk of radius $r$. 

\begin{lemma}[Normalization]
\label{normalization lemma} We have the following properties:
\begin{enumerate}
\item \label{item:normalization1} After a deformation of the family of symplectic fibrations $\pi_t$ produced by \Cref{prp:fibrations}, we may assume that any fiber $\pi_t^{-1}(e^{i\theta})$ above the equator $e^{i\theta} \in S^1 \subset S^2$ coincides with the standard fiber $\op{pr}_2^{-1}(e^{i\theta})$ near $D_\infty \cup L_0$.

\item \label{item:normalization2} After a Hamiltonian isotopy of $L_1$ confined to $S^2 \times S^2 - (D_\infty \cup L_0)$ we may assume that
$$ L_1 \subset S^2 \times S^2 - (D_\infty \cup \ell_1 \cup L_0).$$ 
Moreover, we may assume that the fibrations $\pi_t$ produced by \Cref{prp:fibrations} satisfy $\pi_t=\op{pr}_2$ in some neighbourhood of $\ell_1$ for all $t \in [0,1]$.
\end{enumerate}
\end{lemma}
\begin{proof}

\eqref{item:normalization1}: Recall that \Cref{prp:fibrations} \eqref{item:over-equator} implies that each line $\pi_t^{-1}(e^{i\theta})$ which is the fiber above a point $e^{i\theta} \in S^1 \in S^2$ on the equator coincides with some standard line $\op{pr}_2^{-1}(e^{if_t(\theta)})$ in a small neighbourhood of $L_0$. Consider a one-parameter family of automorphisms of $L_0 = \bb{T}^2$ that connect $f_t \times \op{id}_{S^1}$ to $\op{id}_{\TT^2}$. This family extends to a Hamiltonian isotopy of the Weinstein neighborhood associated to the embedding $\phi_0$ (see \ref{item:disjoint-image}). This isotopy can be further extended by identity to all of $S^2 \times S^2$ after applying a suitable cut-off. By applying this isotopy, we can achieve that $\pi_t^{-1}(e^{i\theta})$ coincides with the standard line $\pi_0^{-1}(e^{i\theta})=\op{pr}_2^{-1}(e^{i\theta})$ near $L_0$ for all $t \in [0,1]$ and $\theta \in S^1$.

It remains to normalize the lines $\pi_t^{-1}(e^{i\theta})$ near $D_\infty$. We begin with the observation that each such line is holomorphic near $D_\infty$ and intersects $D_{\infty}$ transversely in precisely two points. Moreover, it follows from positivity of intersection and automatic transversality (see \cite[Lem.\ 3.3.3]{mcduff-sal-jcurves}) that the intersection of the family of lines $\{\pi_t^{-1}(e^{i \t}) \}_{\t \in [0, 2\pi)}$ with each of the two lines $\{i\} \times S^2$, $i\in \{0,\infty\}$ consists of embedded closed curves $\{i\} \times \gamma^i_t \subset \{i\} \times (S^2 \setminus \{0,\infty\})$ which moreover are homotopically essential. The deformation is performed in two steps:
\begin{enumerate}[label=(\Roman*)]
\item We straighten the lines $\pi_t^{-1}(e^{i\theta})$ so that they coincide with the standard lines $S^2 \times \{\gamma^i_t(\theta)\}$ near $\{i\} \times S^2$, for $i \in \{0,\infty\}$.
\item We choose a family of isotopies of curves $\{\gamma^i_{s;t }\}_{s, t \in [0,1]}$ in $S^2 \setminus \{0,\infty\}$ such that $\gamma^i_{0; t}(\theta)=\gamma^i_t(\theta)$ and $\gamma^i_{1;t}(\theta)=e^{i\theta} \sub S^2$. We use this family to deform the lines $\pi_t^{-1}(e^{i\theta})$ through symplectic lines so that they coincide with $S^2 \times \{\gamma^i_{1;t}(\theta)\}$ near $D_\infty$.
\end{enumerate}

\textit{Step (I):} Fix a symplectomorphism $$\iota_0^i: \{i\} \times (S^2 \setminus \{0,\infty\}) \subset S^2 \times S^2 \to \S \sub (\R^2, \o_0),$$ where $\S$ is an annular domain and $\o_0$ is the standard symplectic form. For $\e$ small enough,  the map $\iota_0^i$ extends to a symplectomorphism $$\iota^i: \mc{U}^i \to (D_{\e} \tms \S, \o_0 \oplus \o_0) \sub (\R^4, \o_0),$$ where $\mc{U}^i \sub S^2 \tms S^2$ is a neighborhood of $\{i\} \times (S^2 \setminus \{0,\infty\})$. We may assume that $\iota^i$ preserves the product structures and that $\op{Im} \iota^0 \cap \op{Im} \iota^{\infty} = \emptyset$. It will be convenient to identify $\mc{U}^i$ with its image under $\iota^i$. 

After a smooth, $C^1$-small perturbation, we may assume that each line $\pi^{-1}_t(e^{i\theta})$ is linear near $\{i\} \times S^2$ with respect to the standard coordinates on $D_{\e} \tms \Sigma \subset \R^4$. (Use the fact that two embeddings which are sufficiently $C^1$-close are isotopic, and then perform an interpolation.) The lines can be assumed to remain symplectic under this perturbation.

The lines $\pi^{-1}_t(e^{i\theta})$ trace out a family $T^i_{t}: \t \mapsto T_{\gamma^i_t(\t)} (\pi^{-1}(e^{i \t})) \in \op{Gr}_2(T^i_{\gamma^i_t(\t)}(\R^4))$ of symplectic $2$-planes over $\gamma^i_t$, where we are viewing $\gamma^i_t$ as a path in $\{0 \} \tms \S $ via $\iota^i$. For $s \in [0,1]$, let $T^i_{s,t}$ be a family of symplectic $2$-planes over $\gamma^i_t$ such that $T^i_{0,t}= T^i_t$ and $T^i_{1,t}$ is orthogonal to $\{0\} \times \S$ with respect to the standard Euclidean metric. 

For $\e' \ll \e$, there exists a family of embeddings $F^i_{s, t}: D_{\e'} \tms S^1 \to D_{\e} \tms \S$ with the following properties:
\begin{itemize}
\item $F^i_{s,t}(0, \t)= \g_t^i(\t)$ for $\t \in S^1$, 
\item  $F^i_{s,t}(-, \t): D_{\e'}  \to \R^4$ is an affine map, 
\item $T ( F^i_{s,t} (D_{\e'}, \t) ) = T^i_{s,t}(\t)$. 
\end{itemize}

For $N$ large enough, the pairs of embeddings $\{F^i_{n/N, t}, F^i_{(n+1)/N, t} \}_{n=0}^{N-1}$ satisfy the conditions of \Cref{interpolation}. By repeated application of (the parametric version of) \Cref{interpolation}, we therefore obtain the desired deformation.

\textit{Step (II):} Let $d(-, -)$ be the distance on $D_{\e} \sub \R^2$ induced by the Euclidean metric. Fix a (not-necessarily symplectic) compactly-supported isotopy $\psi^i_{s, t} \colon \S \to \S$ such that $\psi^i_{0,t}= \op{Id}$ and $\psi^i_{s,t}(\gamma^i_{0,t}(\theta)=\gamma^i_{s,t}(\theta)$. 
For $\delta>0$, let $\rho:[0, \delta] \to \R$ be a non-increasing function such that $\rho(t)=1$ for $t \in [0, \delta/3]$ and $\rho(t)=0$ for $t \in [2\delta/3, \delta]$. 
For $i \in \{0, \infty \}$, consider the family of diffeomorphisms $\mu^i_{t}$ taking
$$ D_{\delta} \tms \S \ni (z_1,z_2) \mapsto  \mu^i_{t}(z_1, z_2)= (z_1,\psi^i_{\rho(d(z_1,0))}(z_2)).$$ We may view $\mu^i_t$ as a self-diffeomorphism of $\mc{U}^i$ by extending it as the identity. 

Observe that $\mu^i_{t} ( (S^2 \tms \{\op{pt} \} ) \cap \mc{U}^i)$ is a symplectic line for all $t$. We now construct the desired deformation of the lines $\pi_t^{-1}(e^{i \t} )$ by replacing $\pi_t^{-1}(e^{i \t}) \cap \mc{U}^i$ by its image under $\mu_t^i$.  This completes step (II). 

To complete the proof, choose an almost-complex structure $\tilde{J}^t$ which is compatible with $\o \oplus \o$, such that $\tilde{J}^t=i$ near $D_{\infty}$ and such that the lines $\pi_t^{-1}(e^{i \t})$ are $\tilde{J}^t$-holomorphic for all $\t \in S^1$. Gromov's theorem \cite[0.2.A]{gromov} implies that the lines $\pi_t^{-1}(e^{i \t})$ extend to a global foliation by $\tilde{J}^t$-holomorphic spheres. The family of fibrations is now defined as in \Cref{prp:fibrations}. 

\eqref{item:normalization2}: After the above deformations, $\pi_t^{-1}(1)$ is an isotopy of symplectic lines that connect $\ell_1=\pi_0^{-1}(1)=\op{pr}_2^{-1}(1)$ to a symplectic line that is disjoint $L_1$ (by Part (4) of \Cref{prp:fibrations}), where these lines moreover remain fixed in some neighbourhood of $D_\infty \cup L_0$. We can now use \cite[Prop.\ 0.3]{siebert-tian} to produce a Hamiltonian isotopy $\phi_t^{H_t}$ that generates $\pi_t^{-1}(1)=\phi_t^{H_t}(\ell_1)$. The sought Hamiltonian isotopy of $L_1$ can then be taken to be $(\phi_{1-t}^{H_t})^{-1}(L_1)$.

Since $L_1$ is disjoint from $\ell_1$ after the Hamiltonian isotopy, we may deform the path $J^t$ given by \Cref{prp:fibrations} to one for which $J^t=i$ is satisfied near $\ell_1$. For this path of almost-complex structure, the additional claim that $\pi_t=\op{pr}_2$ holds in the same neighbourhood can be assumed.
\end{proof}


Given a normed vector space $(V, \|- \|)$ and a pair of subspaces $V_1, V_2 \sub V$, we define $d(V_1, V_2):= \sup_{\|v \|=1} \| (P_{V_1} - P_{V_2})(v) \|$, where $P_{V_i}: V \to V$ is the projection onto $V_i \sub V$.

Let $\op{Gr}_2^{\o}(\R^4) \sub \op{Gr}_2(\R^4)$ be the subset of $2$-planes in $\R^4$ on which the standard symplectic form is non-degenerate. Observe that $d(-,-)$ makes $\op{Gr}_2(\R^4)$, and hence also $\op{Gr}_2^{\o}(\R^4)$, into a metric space.

\lem[Interpolation] \label{interpolation}
Let $\g: S^1 \to \{0 \} \tms \R^2 \sub (\R^4, dx_1\wedge dy_1 + dx_2 \wedge dy_2)$ be an embedded loop. Fix $0 < \e < 1/100$.

Let $\t \mapsto \L_{\t}  \in \op{Gr}_2(\R^4)$ be a path of $2$-planes such that $\L_{\t} \pitchfork \dot{\g}(\t)$ and $B_{100\e}(\L_{\t}) \sub \op{Gr}_2^{\o}(\R^4)$.

Let $S_0, S_1:  D_1 \tms S^1 \to \R^4$ be embeddings of solid tori with the following properties:
\begin{enumerate}
\item \label{item:core} $S_i(0, \t)= \g(\t),$ 
\item \label{item:linear} $S_i(-, \t): D_1 \sub \R^2 \to \R^4$ is a linear map, 
\item \label{item:close} For each $\t \in S^1$, the planes $T(S_i(D_1, \t))$ are contained in the ball $B_{\e}(\L_{\t}) \sub \op{Gr}_2^{\o}(\R^4)$. 
\end{enumerate}

Then there exists a family of embeddings $\tilde{S}_s: D_1 \tms S^1 \to \R^4$, for $s \in [0,1]$, such that for $\e' \ll \e$ small enough: 
\begin{itemize} 
\item[(i)] $\tilde{S}_{s}(0, \t)= \g(\t)$,
\item[(ii)] $ \tilde{S}_{s}((D_1-D_{1/2}), \t)$ is contained in $ S_0(D_1, \t)$ 
\item[(iii)] $\tilde{S}_{1}(D_{\e'},\t)$ agrees with $S_1(D_{\e'},\t)$ 
\item[(iv)] $T_{(x,y, \t)} (\tilde{S}_s(D_1,\t))$ is contained in $B_{50\e}(\L_{\t})$ for all $s \in [0,1]$ and any $(x,y) \in D_1$. In particular, $\tilde{S}_s(D_1,\t)$ is a symplectic submanifold. 
\end{itemize}

Given a closed interval $I \sub \R$, if $\g, S_0, S_1$ depend smoothly on a parameter $\s \in I$, then $\tilde{S}_s$ can also be taken to depend smoothly on $\s \in I$. 

\elem

\pf
Fix an embedding $\phi: D_1 \tms S^1 \tms (-\eta, \eta) \to \R^4$, for some $\eta>0$, such that $\phi(-, -,0)= S_1(-,-)$. Let $\mc{U}= (\op{Im} \phi)^{\circ}$ and observe that $\phi^{-1}$ defines natural coordinates $(x,y, \t, z)$ on $\mc{U}$, where $(x,y) \in D_1$.

For $i=0,1$, consider the solid torus $\mc{T}_i= \op{Im}S_i$, and observe that $\mc{T}_i$ admits a foliation $\mc{F}_i$ with leaves $\{L_i(\t)=  S_i(D_1,\t) \}_{\t \in S^1}$. 

According to \eqref{item:close}, we can parametrize $L_0(\t)$ in terms of the coordinate chart $\phi^{-1}$ in some neighborhood of $\g$. More precisely, for some $\eta_2>0$, we can write:
\eq L_0(\t)=  (x,y, a(x,y,\t), b(\t, x,y)),\eeq
where $a,b: D_{\eta_2} \tms S^1 \to S^1$ are of the form $a(x, y, \t)= \t+ a_{\t}(x,y)$, $b(x,y,\t)= b_{\t}(x,y)$, with $a_{\t}(0,0)= b_{\t}(0,0)=0$.  Observing that $\d_{\t}a_{\t}(0,0)=0$, we can assume after possibly shrinking $\eta_2$ that $|\d_{\t}a_{\t}(x,y)|<1/100$ for all $(x,y) \in D_{\eta_2}$. 

Given $0<\delta<1/100$, let $\rho: [0, \delta] \to \R$ be a non-increasing function such that $| \rho' | < 2/\delta$ and 
\eqs  \rho= \begin{cases} 1 &\text{in } [0, \delta^2],  \\
0 &\text{in } [\delta, 1].
\end{cases} \eeqs

For $s \in [0,1]$, let $\tilde{\rho}_s: D_{\eta_2} \to \R$ be defined as $\tilde{\rho}_s(x,y)=  s \rho ( \|(x, y) \| )$. 

For $\t \in S^1$, consider the map $f_s^{\t}: D_{\eta_2} \to \R^4$ defined by:
$$(x,y) \mapsto (\tilde{\rho}_s a(x,y,\t) + (1-\tilde{\rho}_s) \t, x, y, \tilde{\rho}_s b(x,y,\t)),$$ where the addition in the first factor comes from the group structure on $S^1$. 

Let $L'_s(\t)= \op{Im} f_s^{\t}$. 



Observe that $L_s'(\t)$ agrees with $L_1(\t)$ outside a neighborhood $\op{Op}_{\eta_3}(\g)$, for some $\eta_3>0$. So we let 
\eq \tilde{L}_s(\t)= \begin{cases} L_s'(\t) &\text{in } \op{Op}_{\eta_3}(\g) \\
L_1(\t) &\text{else.} 
\end{cases}
\eeq

\lem
For $\delta>0$ small enough, $\tilde{L}_s(\t)$ is an embedded symplectic submanifold. In fact, $T_p\tilde{L}_s(\t) \sub B_{50\e}(\Lambda_{\t})$ for all $p \in \tilde{L}_s(\t)$. 
\elem
\pf It's enough to show that $L_s'(\t)= \op{Im} f_s^{\t}$ is an embedded symplectic submanifold. It is clear that $\op{Im} f_s^{\t}$ is embedded. It is also clear that $\op{Im} f_s^{\t}$ is symplectic for $\| (x,y) \| \geq \delta$, with $T\op{Im} f_s^{\t}\sub B_{50\e}(\Lambda_{\t})$. 

It remains to show that $T \op{Im} f_s^{\t} \sub B_{50\e}(\Lambda_{\t})$ for $\| (x,y) \| \leq \delta$, for $\delta$ small enough. Let $F_s^{\t}(x,y):= \tilde{\rho}_s(x,y) a(x,y,\t) + (1-\tilde{\rho}_s(x,y)) \t$ and let $G_s^{\t}(x,y)= \tilde{\rho}_s b(x,y,\t)$. By \eqref{item:close}, it is enough to show that $| \d_xF_s^{\t} | < 6\e$ and $|\d_y F_s^{\t} | < 6\e$ and similarly for $G_s^{\t}$.

Observe that we have $a_{\t}(x,y)= a_xx+a_y y + O(\|(x,y) \|)$.  Now given $\|(x_0, y_0) \|< \delta$, we have 
\begin{align*}
| \d_x F_s^{\t} (x_0, y_0) | &=  | s \rho' \d_x(\sqrt{x^2+y^2}) (a(\t, x_0, y_0) - \t) + \tilde{\rho}_s \d_x a(x_0, y_0,\t) | \\ 
&\leq  | 2/ \delta |  ( |a_x| |\delta| + |a_y| |\delta| +  C |\delta^2| ) + |\rho| ( |a_x| + C |\delta^2| ),
\end{align*} 
for some $C>0$ independent of $\delta$. 

According to \eqref{item:close}, we have $|a_x|<\e$ and $|a_y|<\e$. Choosing $\delta$ small enough, we conclude that $| \d_x F_s^{\t} (x_0, y_0) |  \leq 6 \e$. The same argument works for $\d_yF_s^{\t}$ and for $\d_xG_s^{\t}, \d_y G_s^{\t}$. 
\epf

We now define $\tilde{S}_s: D_1 \tms S^1 \to \R^4$ by first setting $\tilde{S}_s(x, y, \t)= f_s^{\t}(x,y)$ for $(x,y) \in D_{\eta_2}$, and by then extending this map so that $\tilde{S}_s(-, \t): D_1 \to \tilde{L}_s(\t)$ is a family of embeddings depending smoothly on $s$ and $\t$. It is clear that this extension exists by definition of $f_s^{\t}$. 

To check that $\tilde{S}_s$ is a smooth embedding, we just need to verify that it is injective with non-singular differential. It's enough to check this for $(x,y) \in D_{\eta_2}$. Using now our above assumption that $|\d_{\t}(x,y)|<1/100$ for $(x,y) \in D_{\eta_2}$, we find that $\d_\t \tilde{S}_s>1/2$ and that $\{\d_\t \tilde{S}_s, \d_x \tilde{S}_s, \d_y \tilde{S}_s \}$ are linearly independent. This proves the claim.
\epf

\begin{lemma}[Inflation]
\label{lemma:liouville-inflation}
There exists a smooth Liouville flow
\begin{gather*}
\psi^s \colon (S^2 \times S^2 - (D_\infty \cup \ell_1)) \hookrightarrow (S^2 \times S^2 - (D_\infty \cup \ell_1)),\\
(\psi^s)^* (\omega \oplus \omega)=e^s(\omega \oplus \omega), \:\:\: \psi^0=\op{id},
\end{gather*}
which is defined for all $s \in (-\infty,0]$ and which satisfies the following properties:
\begin{enumerate}
\item \label{item:inflation1} $\psi^s(L_0 - \ell_1) \subset L_0 - \ell_1$.
\item \label{item:inflation2} Every fiber of $\op{pr}_1$ (resp.\ $\op{pr}_2$) get mapped by $\psi^s$ into another fiber of $\op{pr}_1$ (resp.\ of $\op{pr}_2$).
\item \label{item:inflation3} There exists a compact subset $C \subset S^2 \times S^2 - (D_\infty \cup \ell_1)$ for which  $\psi^s(S^2 \times S^2 - (D_\infty \cup \ell_1)) \subset C$ is satisfied whenever $s \ll 0$ is sufficiently small.
\end{enumerate}
\end{lemma}

\begin{figure}[htp]
\vspace{3mm}
	\labellist
	\pinlabel $m$ at 165 -4
	\pinlabel $m$ at 69 -4
	\pinlabel $L_0$ at 157 15
	\pinlabel $L_0$ at 60 34
	\pinlabel $s$ at 165 30
	\pinlabel $s_1$ at 49 30
	\pinlabel $s_2$ at 90 30
	\pinlabel $\op{pr}_1^{-1}(1)-(D_\infty\cup\ell_1)$ at 70 55
	\pinlabel $\op{pr}_2^{-1}(1)-(D_\infty\cup\ell_1)$ at 165 55
	\endlabellist
	\includegraphics[scale=1.5]{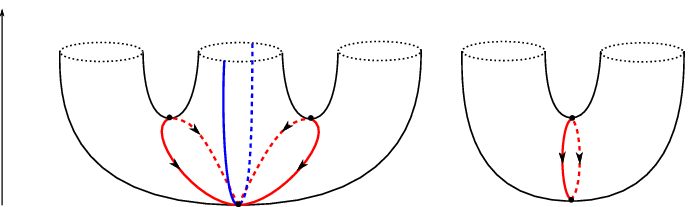}
	\caption{The Liouville flow $\psi^s$ can be taken to be the negative gradient flows on the three and two punctures surfaces induced by the height function showed above. The skeleton, which is fixed by the flow, is shown in red. In addition, the blue curve on the left can also be assumed to be fixed setwise by the flow. The product of the blue curve on the left and the red circle on the right is identified with $L_0 - \ell_1$.}
	\label{fig:liouville}
\end{figure}

\begin{proof}
Observe that the subset
$$ (S^2 \times S^2 - (D_\infty \cup \ell_1),\omega \oplus \omega) \sub (S^2 \tms S^2, \o \oplus \o)$$
is naturally the product of a twice-punctured sphere with a thrice-punctured sphere. 

The desired Liouville flow can be constructed as the product of Liouville flows on the punctured spheres. On the twice-punctured sphere, we can take the pushforward via \eqref{equation:weinstein-s1} of the standard Liouville flow on $D_{1/2}T^*S^1$, which preserves the zero section. 

On the thrice-punctured sphere, a suitable Liouville flow can be constructed as follows. Start with the standard Liouville flow on $(D^2,\omega=\frac{1}{2}d(xdy-ydx))$ which fixes the imaginary part setwise (the imaginary part will be identified with the blue curve shown on the left in Figure \ref{fig:liouville}). Then attach two standard Weinstein one-handles on each of the two embeddings of $S^0 \subset \partial D^2$ contained in the subsets $\{x > 0\} \cap D^2$ and $\{x < 0\} \cap D^2$ respectively. (A standard Weinstein handle in this dimension is just a strip with a suitable Liouville form that agrees with the one on $D^2$ near two of the sides, and which has one non-degenerate critical point of saddle type.) This produces the sought Liouville form, at least up to a suitable symplectomorphism.
\end{proof}

\begin{lemma}
\label{inflated action}
The Liouville flow from \Cref{lemma:liouville-inflation} produces a Lagrangian isotopy
$$ L^s_1 \coloneqq \psi^{s}(L_1) \subset S^2 \times S^2 - (D_\infty \cup \ell_1 \cup L_0) \subset T^*\TT^2, \:\:\: t \le 0,$$
where the symplectic action of $L^s_1$ is of the form
$$\mathbf{p}^{L^s_1}=(0,a_s)$$
for some path $a_s>0$ that depends smoothly on $s \le 0$.
\end{lemma}
\begin{proof}
Since $a_0=a>0$ and since any exact Lagrangian torus in $T^* \bb{T}^2$ intersects the zero section, it is enough to prove the a priori weaker statement that $p^{L^s_1}_1 \equiv 0$, i.e.~that $\mathbf{p}^{L^s_1}=(0,a_s)$ for some arbitrary path $a_s \in \R$. 

The claim $p_1^{L^s_1} \equiv 0$ can be seen to be a consequence of the following fact: there exists a smooth two-chain $C$ inside $S^2 \times S^2 - (D_\infty \cup \ell_1 \cup L_0)$ with boundary in the homology class $\mathbf{e}_1^{L_1}-\mathbf{e}_1^{L_0} \in H_1(L \cup L_0)$. One can readily construct $C$ by appealing to \Cref{prp:fibrations} together with \Cref{normalization lemma}. For example, $C$ can be taken to be a smooth cylinder which lives over an embedded path $\gamma \subset S^2 - \{1\}$ of the base of the fibration $\pi_1$, where $\gamma$ has one boundary point on the embedded circle $\pi_1(L_1)$ and the other boundary point on the equator $S^1=\pi_1(L_0)$.

Given the existence of $C$, the calculation of the symplectic action of $L^s_1$ can now be done as follows. Let $\lambda=p\,dq+\eta$ be the Liouville form on $S^2 \times S^2 - (D_\infty \cup \ell_1)$ that defines the Liouville flow $\psi^s$, where $\eta$ is some closed one-form on $S^2 \times S^2 - (D_\infty \cup \ell_1)$. Observe that $(\psi^s)^*\lambda=e^s\lambda$. Since $L_0$ is a Lagrangian that remains fixed setwise under $\psi^s$ by assumption, one can deduce that $\lambda|_{T(L_0 - \ell_1)} \equiv 0$. Using these facts we finally compute
$$p_1^{L^s_1}=p_1^{L^s_1}-p_1^{L_0} =\int_{\mathbf{e}_1^{L^s_1}-\mathbf{e}_1^{L_0}}\lambda =\int_{\psi^s(C)} d\lambda=e^s \int_C d\lambda=e^s(p_1^{L_1}-p_1^{L_0})=0$$
for all $s \le 0$.
\end{proof}

\subsection{Proof of Proposition \ref{prp:IntoConvex}} \label{subsection:finalproof}
Fix a family $\pi_t: S^2 \tms S^2 \to S^2, t \in [0,1],$ of symplectic $S^2$-fibrations having the properties furnished by \Cref{prp:fibrations} and \Cref{normalization lemma}. Let us consider the smooth one-parameter family of hypersurfaces 
$$Y_t \coloneqq \pi_t^{-1}(S^1) \subset S^2 \times S^2.$$
These hypersurfaces are foliated by the symplectic lines $\pi_t^{-1}(e^{i \theta}),$ $e^{i\theta} \in S^1 \sub S^2$ (recall from \Cref{subsection:notation-reductions} that we have fixed subsets $S^1 \sub S^2, \{0\} \sub S^2, \{\infty\} \sub S^2$ which are respectively called the equator, south pole and north pole). 


The path $\R / \Z \ni \t \mapsto e^{2\pi i \t}$ lifts under $\pi_t$ to the characteristic foliation $\ker (\omega|_{TY_t}) \subset TY_t$. By integrating, we obtain \emph{symplectic monodromy maps} $$\varphi_t: (S^2, \o) \to (S^2, \o),$$ where we identify $(S^2 \tms \{1 \}, \o \oplus \o)= (S^2, \o)$ via $\op{pr}_1$. We may view the monodromy maps as a smooth path $t \mapsto \varphi_t$, where $\varphi_0= \op{id}_{S^2}$ in view of \Cref{prp:fibrations} \eqref{item:t=0}. It follows from \Cref{normalization lemma} \eqref{item:normalization1} that the characteristic foliation is standard near $D_\infty \cup L_0$. This implies that a neighbourhood of $\{0,\infty, S^1\} \sub S^2$ is fixed pointwise by $\varphi_t$.  We thus have that $\varphi_t=\varphi_t^{H_t}$ for some Hamiltonian $H_t \colon S^2=S^2 \times \{1\} \to \R$ which can be taken to be locally constant in a neighbourhood of $S^1 \cup \{0,\infty\}$. Without loss of generality, we may assume that $H_t$ vanishes near $S^1$.

Let $\s<0$ be a negative constant which will be fixed later. For $s \in (-\infty, 0]$, let $\psi_s$ be the Liouville flow on $W:=(S^2 \tms S^2 - (D_{\infty} \cup \ell_1))$ furnished by \Cref{lemma:liouville-inflation}. Let $\l$ be the associated Liouville form. Let $\widehat{W}$ be the canonical completion of the Liouville manifold $(W, \l)$. 

Let us now construct a new family $Y_t^{\s}$ of hypersurfaces by deforming $Y_t$ under the Liouville flow. 

Let $\widehat{Y}_t \sub \widehat{W}$ be the completion of $Y_t \sub W$ (i.e. $\widehat{Y}_t$ is constructed by extending $Y_t$ via the Liouville flow on $\widehat{W}$). Let $\check{Y}_t^{\s}= \psi_{\s}(\widehat{Y}_t) \cap W$, where $\psi: \R \tms \widehat{W} \to \widehat{W}$ is the Liouville flow extended to the completion. It follows from \Cref{normalization lemma} \eqref{item:normalization1} that the $\check{Y}_t^{\\s}$ coincides near $D_{\infty}$ with fibers of the projection $\op{pr}_2$. Similarly, it follows by \Cref{normalization lemma} and \Cref{lemma:liouville-inflation} \eqref{item:inflation2} that $\check{Y}_t^s$ coincides near $\ell_1$ with the fibers of $\op{pr}_2$. We now define $Y_t^{\s}$ as the closure of $\check{Y}_t^{\s} \sub S^2 \tms S^2$.


\rmk For a different perspective on the construction of $Y_t^{\s}$, observe that $\pi_t: \check{Y}_t^{\s} \to S^1-\{1\}$ is a fibration by annuli. We are obtain the $S^2$-fibration $Y_t^{\s} \to S^1$ by gluing in a pair of points to each fiber over $e^{i \t} \in S^1 - \{1\}$, and by gluing in an $S^2$ fiber over $\{1\}$. \ermk

We now study the monodromy map $\varphi_t^{\s}: (S^2, \o) \to (S^2, \o)$ induced by the characteristic flow along $Y_t^{\s}$, where we continue to identify $(S^2 \tms \{1 \}, \o \oplus \o)= (S^2, \o)$ via $\op{pr}_1$. The following lemma shows that the family $t \mapsto \varphi_t^{\s}$ can be generated by a family of Hamiltonians $H^{\s}_t$ whose size depends on $\s$ and on the Hamiltonians $H_t$ generating the monodromy of $Y_t$.

\begin{lemma} \label{lemma:hamiltonian-sup-bound}
The family of monodromy maps $\varphi_t^{\s}, t \in [0,1]$ is generated by a family of Hamiltonians $H^{\s}_t$ which satisfies the bound
$$\max_{t,x} |H^{\s}_t(x)| \le e^{\s}\max_{t,x}|H_t(x)|.$$
\end{lemma}

\pf

According to \Cref{normalization lemma} \eqref{item:normalization2} and the definition of $Y_t^{\s}$, there exists $\e>0$ small enough so that the map $\op{pr}_2^{-1}(e^{-i\e}) \to \op{pr}_2^{-1}(e^{i \e})$ induced by the  characteristic flow along $Y_t^{\s}$ is the identity (i.e. it intertwines the projection $\op{pr}_1$).  Hence it is enough to study the map $\op{pr}_2^{-1}(e^{i \e}) \to \op{pr}_2^{-1}(e^{-i\e})$ induced by the characteristic flow. According to \Cref{normalization lemma} \eqref{item:normalization1}, this map is the identity near $\{0, \infty \} \sub S^2 = \op{pr}_2^{-1}(e^{i \e})$. Hence, we may restrict our attention to the map of annuli $\tilde{\varphi}_t^{\s}: \op{pr}_2^{-1}(e^{i \e}) - D_{\infty} \to \op{pr}_2^{-1}(e^{-i\e}) - D_{\infty}$ induced by the characteristic flow along $Y_t^{\s}- D_{\infty}$.

We have a canonical identification $\op{pr}_2^{-1}(e^{i \e}) - D_{\infty}= S^2 - \{0, \infty \}$, so we can view $\tilde{\varphi}_t$ as a symplectomorphism of $(S^2- \{0, \infty \}, \o)$.

Let us also consider the monodromy of annuli $\tilde{\varphi}_t: (S^2- \{0, \infty\}) \to (S^2- \{0, \infty\})$ induced by the characteristic flow along $Y_t$. We let $\tilde{H}_t$ be the generating Hamiltonian. 

Observe that $\tilde{\varphi}_t^{\s}=\psi^{\s} \circ \tilde{\varphi}_t \circ (\psi^{\s})^{-1}$, where $\varphi_t$ is the Liouville flow on $S^2 \tms S^2 - (D_{\infty} \cup \ell_1)$. Hence $\tilde{\varphi}_t^{\s}$ is generated by a Hamiltonian $\tilde{H}^{\s}_t=e^{\s}H_t \circ (\psi^{\s})^{-1}$. 

The monodromy $\varphi_t^{\s}$ is now generated by a Hamiltonian $H^{\s}_t$ defined by extending $\tilde{H}^{\s}_t: (S^2- \{ 0, \infty \}) \to \R$ to be locally constant near the poles.  
\epf



Let us now analyze the structure of the $S^2$-fibration $\pi_t: Y_t^{\s} \to S^1$. Recall that $\pi^{-1}(e^{i \theta})= \op{pr}_2^{-1}(e^{i \theta})$ if $\t \in [-2\e, 2\e]$, for $\e$ small enough. On the other hand, \Cref{lemma:liouville-inflation} \eqref{item:inflation3} implies that $\psi_{\s} (Y_t - \ell_1)$ and the Lagrangian $L_1^{\s}$ are both contained inside a fixed compact subset $C \subset S^2 \times S^2 - (D_\infty \cup \ell_1)$, provided that $\s$ is small enough. 

Hence, after possibly making $\e$ smaller, we have
\eq \label{equation:backwards-flow-disjoint} \op{Op}_{\e}( \{\pi^{-1}_t(e^{i \theta}) \mid \t \in [-\e/2, \e/2] \} ) \cap (Y_t^{\s} \cup L_1^{\s}) \sub \{ \pi^{-1}_t(e^{i \theta}) \mid \t \in [-2\e, 2\e] \}. \eeq We emphasize that \eqref{equation:backwards-flow-disjoint} holds for any sufficently small choice of $\s$; in particular, $\e$ does not depend on $\s$.

Given $\delta$ small enough (depending on $\e$), we have a symplectic embedding
\eq \label{equation:suspension-embedding} (S^2 \tms (-\delta, \delta)^2, \o \oplus dx \wedge dy) \to \op{Op}_{\e}( \{\pi^{-1}(e^{i \theta}) \mid \t \in [-\e/2, \e/2] \} ) \sub S^2 \tms S^2,\eeq which is the identity on the first factor and is defined by \eqref{equation:weinstein-s1} in the second factor. Observe that $S^2 \tms \{y=0\}$ embeds into $Y_t^{\s}$ under \eqref{equation:suspension-embedding}.

Let $\beta \colon [-\delta, \delta] \to [0,1]$ be a smooth bump function which satisfies $\beta' \ge 0$, and $\beta(t)=0$ near $\{t=-\delta\}$ while $\beta(t)=1$ near $\{t=\delta \}$.

For $\tau \in [-\delta, \delta]$ and $t \in [0,1]$, let $$\varphi_{\tau;t}^{\s}:= (\varphi_{t \beta(\tau)}^{\s})^{-1}: S^2 \to S^2.$$ 

\begin{lemma}
\label{characteristic action}
For fixed $t \in [0,1]$, the path $\tau \mapsto \varphi_{\tau;t}^{\s}$ of Hamiltonian diffeomorphisms can be generated by the Hamiltonian
$$ K_{\tau; t} \coloneqq -t\beta'(\tau)\cdot H_{t\beta(\tau)} \circ \phi_{t\beta(\tau)}^{H} \colon S^2 \to \R$$ 
\end{lemma}
\begin{proof}
This follows from the well-known formula for the inverse path of a Hamiltonian isotopy, together with the fact that rescaling time has the effect of an analogous rescaling of the Hamiltonian.
\end{proof}

It follows from \Cref{lemma:hamiltonian-sup-bound} and \Cref{characteristic action} that we may choose $\s$ small enough so that $\op{max} | K_{\tau; t}| < \delta$. We view $\s$ as fixed the remainder of this section.

We will now construct a new family hypersurfaces $\tilde{Y}_t$ by applying the so-called \emph{suspension construction}.  Let us first identify $\{(\zeta, x, 0) \mid \zeta \in S^2, x \in (-\delta, \delta) \}$ with $\pi_t^{-1} ( \{ e^{i \t} \mid \t \in I \sub (-\e, \e) \}) \sub Y_t^{\s}$ via \eqref{equation:suspension-embedding}. We let $\mathring{Y}_t^{\s}:= Y_t^{\s} - \{(\zeta, x, 0) \mid \zeta \in S^2, x \in (-\delta, \delta) \}.$ 

We now define 
\eq \tilde{Y}_t= \mathring{Y}_t^{\s} \cup \{( \varphi_{x; t}(\zeta), x, K_{x; t}(\zeta)) \mid \zeta \in S^2, x \in (-\delta, \delta) \}. \eeq

It's straightforward to calculate (cf.\ \cite[Lem.\ 6.7]{dgi}) that the symplectic monodromy map associated to the fibration $\pi_t: \tilde{Y}_t \to S^1$ is the identity.

\lem \label{lemma:hamiltonian-vanishing}
The Hamiltonian $H_t$ vanishes near $S^1 \cup \{0\} \cup \{\infty\}$ for $t \in [0,1]$ (and hence so does $K_{\tau;t}$). 
\elem
\pf Recall that $H_t$ is locally constant near $S^1 \cup \{0\} \cup \{\infty\}$ and was assumed to vanish near the equator $S^1$. It thus suffices to prove that $H_t$ also vanishes near the poles. We'll prove that $H_t$ vanishes near $\{0\}$ since the other pole can be handled analogously. 

By hypothesis, $H_t$ is constant is some neighborhood $\mc{U}$ of $\{0\}$. 

Let $\tilde{Y}_t^{\circ}:= \tilde{Y}_t  \cap (S^2 \tms S^2 - D_{\infty})$. Note that $(S^2 \tms S^2 - D_{\infty})$ is an exact symplectic manifold via the exact symplectomorphism \eqref{equation:weinstein-t2}. Observe that $\tilde{Y}_t$, and hence also $\tilde{Y}_t^{\circ}$ (since $\tilde{Y}_t \cap D_{\infty}$ is characteristic), is foliated by closed characteristic curves. The closed characteristics of $\tilde{Y}_t^{\circ}$ must all have the same action. Since the closed characteristics contained in $L_0= S^1 \tms S^1$ evidently have vanishing action, it follows that all closed characteristics have vanishing action.

Let $\g_t \sub (\mc{U}- \{0\}) \tms S^1 \sub \tilde{Y}_t^{\circ}$ be a closed characteristic. The action of $\g_t$ can be computed as $\int_{-\delta}^{\delta} K_{x;t}(p) dx=-t \int_0^1\tilde{H}_{ts}(p) ds= -\int_0^t \tilde{H}_s(p) ds$, for $p \in \mc{U}- \{0\}$. It follows that $\int_0^t \tilde{H}_s(p) ds=0$, which implies that $\tilde{H}_t$ and hence also $H_t$ vanishes identically in $\mc{U}$. 
\epf

The following corollary following immediately from \Cref{lemma:hamiltonian-vanishing} and the definition of $\tilde{Y}_t$.

\cor
For some $\e>0$ and all $t \in [0,1]$, we have \eq \op{Op}_{\e}(D_{\infty}) \cap \tilde{Y}_t= \op{Op}_{\e}(D_{\infty}) \cap Y_t = \op{Op}_{\e}(D_{\infty}) \cap Y_0. \eeq 
\ecor
\qed

\lem
	There exists a global Hamiltonian $G_t: S^2 \tms S^2 \to \R$, for $t \in [0,1]$, whose flow fixes $D_{\infty}$ setwise and induces the isotopy $\tilde{Y}_t$.
\elem

\pf
We have shown that the characteristic distribution of the hypersurfaces $\tilde{Y}_t$ is the push-forward of a constant vector field for a suitable family of parametrisations. Hence we can extend $\tilde{Y}_t$ to a smooth isotopy of $S^2 \tms S^2$ which preserves the symplectic form in a small neighborhood of $\tilde{Y}_t$ and which is the identity near $\tilde{Y}_t \cap D_{\infty}$. The symplectic action of the closed characteristics of $\tilde{Y}_t$ is independent of $t$. It thus follows by Banyaga's isotopy extension theorem \cite[Thm.\ II.2.1.]{banyaga} that $\tilde{Y}_t$ is generated by a Hamiltonian $G^1_t: S^2 \tms S^2 \to \R$ which is locally constant for fixed $t$ near $\tilde{Y}_t \cap D_{\infty}$. 

We now define $G^2_t: \op{Op}_{\e}( \tilde{Y}_t \cup S^2 \tms \{0 \} \cup S^2 \tms \{ \infty \}) \to \R$ by requiring that $G_t= G^1_t$ in $\op{Op}_{\e}(\tilde{Y}_t)$ and that $G_t$ is constant for each fixed $t$ in $\op{Op}_{\e}(S^2 \tms \{0 \})$ and in $\op{Op}_{\e}(S^2 \tms \{ \infty \})$. 

After possibly shrinking $\e$, we may then define $G^3_t: \op{Op}_{\e}( \tilde{Y}_t \cup D_{\infty}) \to \R$ so that the following properties are satisfied:
\begin{itemize}
	\item[(i)] $G^2_t=G^3_t$ in their common domain of definition,
	\item[(ii)] the Hamiltonian flow associated to $\{G^3_t\}$ is well defined in $\op{Op}_{\e}(\{0 \} \tms S^2 \cup \{\infty\} \tms S^2 )$ and the divisor $\{0 \} \tms S^2 \cup \{\infty\} \tms S^2$ is preserved setwise.
\end{itemize}
It is not difficult to construct a Hamiltonian which satisfies (i) and (ii). Observe however that we cannot in general assume that the divisor $\{0 \} \tms S^2 \cup \{\infty\} \tms S^2$ is preserved \emph{pointwise} by the flow of $G^3_t$ while satisfying (i). 

To complete the proof of the lemma, we simply let $G_t: S^2 \tms S^2 \to \R$ be an arbitrary extension of $G^3_t$. By construction, $G_t$ generates the isotopy $\tilde{Y}_t$ and fixes $D_{\infty}$ setwise.
\epf

We now complete the proof of \Cref{prp:IntoConvex}. By \Cref{inflated action}, we may as well assume that $L= L^{\s}_1$. However, noting that $\tilde{Y}_0=Y_0= \op{pr}_2^{-1}(S^1)$, we find that the Hamiltonian isotopy $(\phi^{1-t}_{G_t})^{-1}(L^{\s}_1)$ takes $L^{\s}_1$ into one of the two components of
$$S^2 \times S^2 - (D_\infty \cup \op{pr}_2^{-1}(S^1)) \subset \{p_2 \neq 0 \}.$$
Each of these components is of the form $\TT^2 \times U$ for a convex $U \subset \R^2- \{0\}$. This proves \Cref{prp:IntoConvex}.
\qed

\section{Stretching the neck along Lagrangian tori} \label{section:neck}

The remainder of this paper is devoted to studying linking of Lagrangians in symplectic $4$-manifolds up to Lagrangian isotopy or smooth isotopy. In this section, we carry out an analysis of pseudoholomorphic foliations of $(S^2 \tms S^2, \o \oplus \o)$ under neck-stretching along multiple disjoint Lagrangian tori. 

This analysis can be viewed as a generalization of \cite[Sec.\ 5]{dgi} from the case of a single Lagrangian torus to the case of multiple disjoint tori. Many arguments are therefore imported directly from \cite[Sec.\ 5]{dgi}. However, the present analysis also differs from \cite[Sec.\ 5]{dgi} in several important respects. The main reason for this is that there are new configurations of broken holomorphic curves which appear in our setting which could be ruled out in \cite[Sec.\ 5]{dgi} for elementary topological reasons.    

Let us now consider a union $L=L_1 \cup \dots \cup L_n \sub  (\R^4, \o)$ of pairwise disjoint Lagrangian tori. Our goal is to prove the following proposition. 

\prop \label{proposition:solidbound}
After possibly relabeling the $L_i$, we can assume that $L_1$ bounds a solid torus in $\R^4- \cup_{j>1} L_j$. Moreover, this solid torus is foliated by symplectic disks. In case the $L_i$ are monotone with monotonicity constant $\kappa_i>0$, we can assume that $\kappa_1 \leq \kappa_j$ for all $j \geq 1$. \eprop

We will also sketch a proof of the following enhancement of \Cref{proposition:solidbound}.

\prop \label{proposition:enhancement} 
After possibly isotoping $L=L_1 \cup\dots\cup L_n$ through Lagrangian tori, we can assume that the solid torus produced by \Cref{proposition:solidbound} has the following additional property: the monodromy map induced by the characteristic flow on the foliation by symplectic disks is the identity. 
\eprop

The proof of \Cref{proposition:solidbound} will occupy the rest of this section and uses both ``hard" and ``soft" tools.

\subsection{A neck-stretching datum} \label{subsection:neckstretch}   Choose $R>0$ sufficiently large so that the polydisk $\mc{P}(R,R) = \{ (z_1, z_2) \in \C^2 \mid |z_1|< R, |z_2|<R\}$ contains $L= L_1 \cup\dots\cup L_n$. We compactify $\mc{P}(R,R)$ to $(S^2 \tms S^2, \o_R \oplus \o_R)$, where $\o_R$ is a rescaling of the standard symplectic form such that $\o_R(S^2)=\pi R^2$. Let $D_{\infty}:= S^2 \tms S^2 - \mc{P}(R,R)$ be the divisor at infinity.

For $i=1,\dots,n$, let $\phi_i: N(0_{\bb{T}^2}) \to N(L_i) \sub S^2 \tms S^2 - D_{\infty}$ be embeddings of Weinstein neighborhoods with disjoint images. Here $N(0_{\bb{T}^2})$ and $N(L_i)$ are fixed open neighborhoods of $0_{\bb{T}^2} \subset T^*\bb{T}^2$ and $L_i \subset \mc{P}(R,R)$ respectively, and $N(L_i) \cap N(L_j) = \emptyset$ if $i \neq j$. We let $\phi= (\phi_1,\dots,\phi_n): N(0_L) \to N(L)$ be the induced Weinstein embedding. 

Let $g$ be a rescaling of the Euclidean metric on $\bb{T}^2= \R^2/\Z^2$.  By choosing this rescaling suitably, we can assume that $S^*_{\leq 4} \bb{T}^2$ is contained in $N(0_{\bb{T}^2})$.  Let $J_{\op{std}}$ and $J_{\op{cyl}}$ be the almost-complex structures on $T^* \bb{T}^2$ and $\R \tms S^*_{1} \bb{T}^2$ constructed in \cite[Sec.\ 4]{dgi}, which were already considered in \Cref{subsection:preliminary-neck-stretch}. Fix a compatible almost-complex structure $J_{\infty}$ on $S^2 \tms S^2 -L$ whose pullback under $\phi_i$ coincides with $J_{\op{cyl}}$ on $S^*_{\leq 4, g} \bb{T}^2-0_{\bb{T}^2}$. 

By choosing $J_{\infty}$ generically in the complement of $N(L)= N(L_1)\cup \dots \cup N(L_n)$ we can assume that it is regular for all simply-covered (possibly punctured) holomorphic curves (see \cite[Thm.\ 7.2]{wendlsft}. Indeed, this follows from the fact that every $J_{\infty}$ holomorphic curve must intersect the complement of $\ov{N(L)}$. 

We now let $\eta= (g, \phi, J_{\op{cyl}}, J_{\op{std}}, J_{\infty}) \in \mc{D}(S^2 \tms S^2, \o_R \oplus \o_R; L)$ be a neck-stretching datum.\footnote{Strictly speaking, to be consistent with \Cref{definition:neck-stretch-datum}, we should view $J_{\op{cyl}}$ and $J_{\op{std}}$ as being defined in the cotangent bundle of the disjoint union $\bb{T}^2 \sqcup\dots\sqcup \bb{T}^2$ ($n$-times).}

\subsection{Index conventions}

With $L \sub S^2 \tms S^2=X$ as in the previous section, we define $c_1^{rel}(-) \in H^2(S^2 \tms S^2, L; \Z)$ as follows. Fix a compatible complex structure on $T^*X$. Given an embedded surface $\S$ with $\d \S \subset L$, choose a generic section $s$ of $\wedge^2_{\C} T^*X|_{\S}$ which belongs to the totally real sub-bundle $\wedge^2(T^*L) \sub \wedge^2_{\C} T^*X$ over $\d \S$, and which moreover is nonvanishing there. Then $c_1^{rel}([\S])$ is the algebraic count of zeros of $s$. 

We define the Maslov class $\mu \in H^2(X,L; \Z)$ by $\mu:= 2c_1^{rel}$.  One can show that in the case where $L$ is connected (so one of the $L_i$ is empty), the pullback of $\mu$ via the Hurewicz map $\pi_2(X, L) \to H_2(X, L; \Z)$ defines a map $\pi_2(X,L) \to \Z$ which agrees with the ``usual" Maslov class; see \cite[Lem.\ 9.2]{cote}.

Finally, if we view $X-L$ as a symplectic manifold with negative cylindrical ends $(-\infty,a] \tms S^*\bb{T}^2$, then we can consider $c_1^{\Phi}(-)$ which is defined as in \cite[Sec.\ 3.1]{dgi}.

\prop[Sec.\ 3.1 in \cite{dgi}] \label{proposition:compactformula} Given a punctured curve $u: \S \to X-L$ with asymptotic orbits $\g_1,\dots,\g_n$, then we have $$2c_1^{\Phi}(u)= c_1^{rel}(\ov{u}) + \sum_{i=1}^n(\mu_{CZ}^{\Phi}(\g_i+\delta)+1),$$ where $\delta>0$ is sufficiently small. \eprop

\subsection{Start of the proof and index analysis} \label{subsection:start-of-stretch-proof}   

Let us now stretch the neck along $L \sub S^2 \tms S^2$ using the datum $\eta \in \mc{D}(S^2 \tms S^2, \o_R \oplus \o_R; L)$, according to the procedure described in \Cref{section:background-material}. We obtain a sequence $J_\tau, \tau \geq 0$, of almost-complex structures on $S^2 \tms S^2$. As usual, Gromov's result \cite[0.2.A]{gromov} ensures that $S^2 \tms S^2$ is foliated by $J_{\tau}$-holomorphic spheres in the class $[S^2 \tms *]$ for all $\tau \geq 0$. We now analyze the SFT-limit of these $J_\tau$-holomorphic spheres as $\tau \to \infty$. 

Throughout this section, we always consider split spheres in the split symplectic manifold $(S^2 \tms S^2 - L) \cup T^*L$ which represent the class $[S^2 \tms *] \in H_2(S^2 \tms S^2; \Z)$. 

Before beginning the main analysis in \Cref{subsection:main-split-sphere-analysis}, we collect some useful properties which are drawn from \cite[Sec.\ 3]{dgi}. 

\prop \label{proposition:fredholm} Suppose that $\bf{u}$ is a split sphere for a regular almost-complex structure $J_\infty$. Then all punctured curves in the top level of $\bf{u}$ have non-negative Fredholm index.  Moreover, all planes in the top level have Fredholm index at least $1$, with equality holding if and only if the plane is simply covered and compactifies to a Maslov $2$ disk. \eprop 

\pf Since $J_{\infty}$ is regular for simply-covered punctured curves, this follows from the proof of \cite[Lem.\ 3.3]{dgi}.  \epf

\prop \label{proposition:punctures} Suppose that $\bf{u}$ is a split sphere for a regular almost-complex structure $J_\infty$. Then the sum of the Fredholm indices of the top-level components of $\bf{u}$ is $2$. Moreover, all components in the building have exactly one or two punctures. All top-level components with two punctures have Fredholm index $0$. \eprop  

\pf The fact that the sum of the Fredholm indices of the top level components is at most $2$ follows from the analysis of \cite[Prop.\ 3.5]{dgi}. For topological reasons, any building must have at least two planes, which must all be in the top level. Hence it follows from \Cref{proposition:fredholm} that the total index of the top level planes is at least $2$. It then also follows from the non-negativity of the index in \Cref{proposition:fredholm} that all other top level components must have Fredholm index zero. 

The fact that all components of the building have exactly one or two punctures also follows from the proof of \cite[Prop.\ 3.5]{dgi}. The basic argument is as follows: if there were a component of the building having three or more punctures, then this would imply (due to the fact that the building has genus $0$) that the building has three or more planes, which is impossible in view of the previous paragraph.   \epf

\cor \label{corollary:planemaslov} Supposing as in the previous proposition that $\bf{u}$ is a spit sphere, then all planes in the top level of $\bf{u}$ compactify to Maslov 2 disks.  In particular, all planes have simply-covered asymptotic orbits.  Moreover, all twice punctured spheres compactify to cylinders which have vanishing Maslov class. \ecor  

\pf The Fredholm index of a punctured sphere $u: \dot{\S} \to S^2 \tms S^2 -L$ satisfies $$\op{ind}(u)= -\chi(\dot{\S})+ \mu(\ov{u}),$$ where $\ov{u}$ is the compactification of $u$; see \Cref{proposition:compactformula}.  Supposing first that $u$ is a plane, it follows by combining \Cref{proposition:fredholm} and \Cref{proposition:punctures} that $\mu(\ov{u})=2$.  If $u$ is a twice punctured sphere, then it follows from \Cref{proposition:punctures} that $u$ has Fredholm index $0$, and hence $\mu(\ov{u})= 0$ as claimed.

Since a plane of Maslov index two that is disjoint from the divisor $D_\infty$ must have a simply covered orbit, we deduce that all remaining components of the building $\bf{u}$ also must have simply covered orbits. (Recall that there are no contractible closed geodesics for the flat metric.) \epf

\subsection{Analysis of the split spheres}  \label{subsection:main-split-sphere-analysis}

We now begin our analysis of the split spheres which arise from stretching the neck along $L=L_1 \cup\dots\cup L_n$. 

\prop \label{proposition:maslovexist} For $i=1,\dots,n$, there exists a dense subset $\mc{U}_i \sub L_i$ and a Maslov $2$ class $\zeta_i \in H_1(L_i; \Z)$ (relative to the natural trivialisation of $T\mc{P}(R,R)$) such that the following property is satisfied: if $p \in \mc{U}_i$, then any split sphere whose compactification passes through $p$ is of Type I, and has all of its asymptotic orbits representing the classes $\pm \zeta_i$. Moreover, such a split sphere exists. 

\eprop 

\rmk The sets $\mc{U}_i$ which we will exhibit are not merely dense: in fact, the $\mc{U}_i$ have full Lebesgue measure and the points of $\mc{U}_i$ can be thought of as generic. More precisely, it will follow from the proof that the $L_i - \mc{U}_i$ is a countable union of geodesics for the flat metric on $L_i$. \ermk

According to \Cref{proposition:maslovexist}, we can introduce the following assumption which will be in force throughout the rest of the proof of \Cref{proposition:solidbound}, i.e. until the end of Section \ref{subsection:inflation}. 

\begin{assumption} \label{assumption:ordering} We assume that the labeling $L= L_1\cup\dots\cup L_n$ is chosen so that $\o(\zeta_i) \leq \o(\zeta_j)$ whenever $i \leq j$ (here, we view $\z_i$ as a class in $\pi_2(\mc{P}(R,R), L_i)$ via the isomorphisms $\pi_2(\mc{P}(R,R), L_i) \to H_2(\mc{P}(R,R), L_i) \to H_1(L_i)$). Observe that if the $L_i$ are monotone as tori in $\mc{P}(R, R) \sub \R^4$, then this is equivalent to the assuming that $\kappa_i \leq \kappa_j$ whenever $i \leq j$, where $\kappa_i$ is the monotonicity constant of $L_i$. \end{assumption} 

We emphasize that this labeling depends on the choice of neck-stretching sequence; it is not intrinsic to the Lagrangians unless they are monotone in $\mc{P}(R,R)$. 

A necessary step towards proving \Cref{proposition:maslovexist} is the following technical lemma. 

\lem \label{lemma:countablegeodesics} For $i=1,\dots,n$, there are at most countably many Reeb orbits of $S^*L_i$ which occur as the asymptotic orbits of a $J_{\infty}$-holomorphic cylinder of Fredholm index zero. \elem

\pf Recall the standard functional analytic setup for punctured holomorphic curves with Morse-Bott asymptotic orbits, as described for instance in \cite[Sec.\ 3.2]{wendlautomatic}. As usual, we write $L= L_1 \cup\dots\cup L_n$. We let $\dot{\S}= \bb{CP}^1- \{0, \infty\}$ and let $W= S^2 \tms S^2-L$.  A homology class $\G \in H_1(L; \Z)$ determines a pair of Morse-Bott manifolds of Reeb orbits $P_{\G}$ of $S^*L$. 

For $p>2$ and $\delta>0$ small enough, we consider the separable Banach manifold $\mc{B}= \mc{B}^{1,p, \delta}(\dot{\S}, W; P_{\G})$. Let $\mc{M}(P_{\G}) \subset \mc{B}$ be the moduli space of simply-covered $J_{\infty}$-holomorphic cylinders. Since $J_{\infty}$ is regular for simply-covered curves, the moduli space $\mc{M}(P_{\G})$ is a smooth manifold whose dimension at a point $u \in \mc{M}(P_{\G})$ is the Fredholm index of $u$. 

Let $\mc{M}(P_{\G})_0 \sub \mc{M}(P_{\G}) \sub \mc{B}$ be the submanifold of curves of Fredholm index zero. Observe that $\mc{M}(P_{\G})_0 \sub \mc{B}$ is a discrete subset. Since $\mc{B}$ is separable, it follows that $\mc{M}(P_{\G})_0$ is countable. It follows that there are at most countably many Reeb orbits occurring as orbits of a cylinder in $\mc{M}(P_{\G})_0$.  Since there is a countable choice of pairs $P_{\G}$, the lemma follows. \epf

\rmk The fact that $\mc{B}$ is separable is not explicitly stated in \cite[Sec.\ 3.2]{wendlautomatic}, but it is not hard to verify. The Banach manifolds which one meets in holomorphic curve theory are usually separable since the Sard-Smale theorem requires separability. \ermk

\pf[Proof of \Cref{proposition:maslovexist}] For $i=1,\dots,n$, let $\mc{U}_i$ be the complement of the union of the geodesics which occur as the projection of asymptotic orbits of a cylinder of Fredholm index zero. It follows from \Cref{lemma:countablegeodesics} that each $\mc{U}_i$ is dense. To see that there exists a split sphere passing through any $p \in \mc{U}_i$, it is enough to observe that there exists a $J_{n}$-holomorphic sphere passing through $p$ for all $n \geq 0$. The limit of such spheres under neck-stretching is the desired split sphere.  

If $\bf{u}$ is a split-sphere whose compactification passes through some $p \in \mc{U}_i$, then it follows from \Cref{proposition:punctures} and the definition of $\mc{U}_i$ that the top-level components of $\bf{u}$ consists of two index $1$ planes and nothing else (in particular, there are no cylinders). These planes must be joined by a bottom-level cylinder, from which we conclude that $\bf{u}$ is of Type I. These planes compactify to Maslov 2 disks according to \Cref{proposition:compactformula}. 

To complete the proof, let us suppose that $\bf{u}_1$ and $\bf{u}_2$ are Type I spheres whose compactification passes through $p_1, p_2 \in \mc{U}_i$, and which have asymptotic orbits representing the Maslov $2$ classes $\eta_{1,i}$ and $\eta_{2,i}$ respectively, where $\eta_{1,i}, \eta_{2,i} \in H_1(L_i; \Z)$. Suppose for contradiction that $\eta_{1,i} \neq \eta_{2,i}$. Observe for $j=1,2$ that $\bf{u}_j$ has a bottom-level cylinder $C_j$ with asymptotic orbits representing the classes $\pm \eta_{j,i}$.  A full classification of such cylinders is described in Section 4 of \cite{dgi}. In particular, according to \cite[Cor.\ 4.3]{dgi}, any two bottom-level cylinders whose asymptotic orbits are not colinear intersect non-trivially in a discrete set.  

This can be seen to give a contradiction by appealing to \cite[Lem.\ 5.8]{dgi}. Indeed, since the Maslov $2$ classes $\eta_{1,i}$ and $\eta_{2,i}$ are distinct, they are not colinear and it follows that the $C_j$ intersect non-trivially in a discrete set.  By positivity of intersection, these intersections are all positive.  Lemma 5.8 of \cite{dgi} now allows us to glue the components of $\bf{u}_1$ and $\bf{u}_2$ to obtain two cycles representing the class $[S^2 \tms *] \in  H_2(S^2 \tms S^2; \Z)$ and intersecting positively. This is the desired contradiction. \epf

\prop \label{proposition:orbitsclass}  For $i=1,\dots,n$, let $\zeta_i \in H_1(L_i;\Z)$ be as above. Then all split spheres in class $[S^2 \tms *]$ have asymptotic orbits which represent $\pm \zeta_i$. \eprop 

\pf Given a split sphere $\bf{u}$, \Cref{proposition:maslovexist} guarantees that we can choose Type I split spheres with orbits in the classes $\pm \zeta_i$ which do not have any orbits in common with $\bf{u}$. We may therefore appeal to \cite[Lem.\ 5.8]{dgi}, from which it follows that all the orbits of $\bf{u}$ represent nonzero multiples of the classes $\pm \zeta_i$. 

It remains to prove that these orbits must in fact represent the classes $\pm \zeta_i$.  It is enough to show that all orbits represent Maslov $\pm 2$ classes. To this end, note that at most one component of the split sphere intersects $D_{\infty}$.

Note that a Maslov $0$ cylinder inside $\R^4 - L$ is asymptotic to two geodesics whose Maslov indices come with different signs (with respect to the trivialisation of $\R^4$). Since the top level components consist of precisely two Maslov $2$ planes together with a number of Maslov $0$ cylinders by \Cref{corollary:planemaslov}, the claim now follows. \epf

\defi A $J_{\infty}$-holomorphic Maslov $2$ plane $u: \C \to S^2 \tms S^2- L$ that is asymptotic to a closed geodesic on $L_1$ in the class $\zeta_1$ and which is disjoint from $D_\infty$ is said to be a \emph{small plane}. We let $\mc{M}_s(J_{\infty})$ denote the moduli space of small planes. \edefi

\lem \label{lemma:cylindergeodesic} If $\g$ is a closed geodesic on $L_i$ which represents the class $\zeta_i$, then there there is a split sphere having a bottom-level cylinder with two positive ends asymptotic to $\pm \g$.  \elem 
\pf We argue as in \cite[Lem.\ 5.12]{dgi}. Fix $p \in \g$. Consider the sequence of $J_l$-holomorphic spheres $u_l$ passing through $p$ and extract a subsequence converging to some building $\bf{u}$. We already argued in \Cref{proposition:orbitsclass} that the asymptotic orbits of $\bf{u}$ must all represent $\pm \zeta_i$ for $i=1,\dots,n$. According to the classification of holomorphic cylinders in \cite[Lem.\ 4.2]{dgi}, the only bottom-level cylinders having orbits colinear to $\zeta_i$ which intersect $p$ must in fact have orbits $\pm \g$. Since all bottom level components are cylinders by \Cref{proposition:punctures} the statement follows. \epf

The following two propositions rely crucially on \Cref{assumption:ordering}. 

\begin{figure}[htp]
\centering
\vspace{3mm}
\labellist
\pinlabel $\color{blue}L_1$ at 0 10
\pinlabel $\color{blue}L_2$ at 160 10
\pinlabel $\color{blue}L_3$ at 40 58
\endlabellist
\includegraphics[scale=1.5]{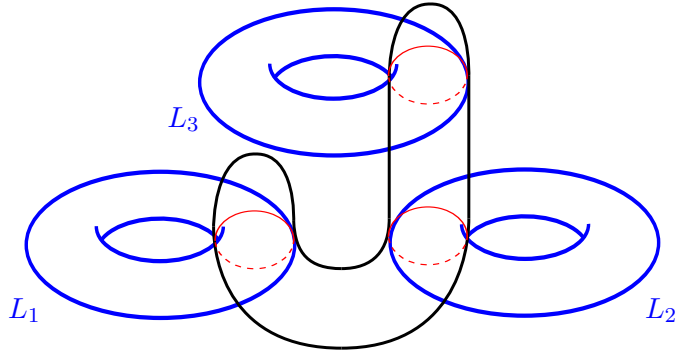}
\caption{A split sphere of Type II passing through $L_1, L_2, L_3$: notice the small plane with asymptotic boundary on $L_1$.}
\label{picture:bigbuilding}
\end{figure}

\prop \label{proposition:smallsplit} Every geodesic of $L_1$  in the class $\zeta_1$ has an embedded small plane asymptotic to it and occurring as a component of a split sphere; see \Cref{picture:bigbuilding}. \eprop  
\pf Let $\g$ be a geodesic in the class $\zeta_1$. According to \Cref{lemma:cylindergeodesic}, there is a bottom-level cylinder $C$ having a positive orbit asymptotic to $\g$. If we remove $C$ from the building, then the remaining components of the building glue together to form a disjoint union of two disks $D_1$ and $D_2$. Only one of these disks can intersect $D_{\infty}$. Without loss of generality, let's assume that $D_1$ doesn't intersect $D_{\infty}$. We claim that $D_1$ doesn't intersect $L_i$ for $i=1,\dots,n$ except at its boundary. Indeed, $D_1$ intersected $L_i$ in an interior circle, then it would have symplectic area strictly greater than $\omega(\zeta_1)$, which is not possible. Hence $D_1$ must be the compactification of a single plane which is asymptotic to $\g$. The embeddedness of this plane follows from positivity of intersection, as in the second paragraph of the proof of \cite[Lem.\ 5.12]{dgi}. \epf

\prop \label{proposition:compactmoduli} The moduli space of small planes is compact. \eprop
\pf First of all, note that the small planes have symplectic area $\o(\zeta_1)$. Hence, they satisfy the appropriate energy bounds for applying the SFT compactness theorem; see the appendix of \cite{cote} for details regarding this standard fact. It follows that any sequence $\{u_i\}_{i =0}^{\infty}$ admits a subsequence converging to a building $\bf{u}$ which is asymptotic to some geodesic $\g$ representing the class $\zeta_1$. We wish to show that $\bf{u}$ is in fact a small plane.

According to \Cref{proposition:smallsplit}, there is a small plane $v$ with orbit asymptotic to $\g$ which occurs as a component of some split sphere $\bf{v}$. Arguing now as in the second paragraph of the proof of Proposition 5.11 in \cite{dgi}, we observe that the union $\bf{u} \cup (\bf{v}- v)$ is a split sphere in the class $[S^2 \tms *]$. Applying now \Cref{proposition:orbitsclass} to the split sphere $\bf{u} \cup (\bf{v}- v)$, it follows that the components of $\bf{u}$ have orbits in the classes $\zeta_i \in H_1(L_i;\Z)$. Since $\o(\zeta_i) \geq \o(\zeta_1)$ by \Cref{assumption:ordering}, we conclude by area considerations that $\bf{u}$ is in fact a small plane. \epf

\subsection{End of the proof of \Cref{proposition:solidbound}} \label{subsection:inflation}

By the compactness of the moduli space of small planes established in \Cref{proposition:compactmoduli}, there is a neighborhood $\mc{N}$ of $L -L_1$ which does not intersect any of the small planes. Let $\tilde{J}_{\infty}$ be an arbitrary extension of $J_{\infty}$ over $\mc{N}$. We now repeat the whole neck-stretching procedure with different data: namely, we replace $J_{\infty}$ with $\tilde{J}_{\infty}$ and we stretch along $L_1$ rather than along $L= L_1 \cup\dots\cup L_n$.  

We are now precisely in the setting of neck-streching along a \emph{single} Lagrangian torus which was considered in \cite[Sec.\ 5]{dgi}. We can therefore appeal directly to the results of their analysis.  Observe first that there is a natural embedding \eq \label{equation:smallplaneembedding} \mc{M}_s(J_{\infty}) \hookrightarrow \mc{M}_s(\tilde{J}_{\infty}). \eeq  Moreover, it is shown in \cite[Lem.\ 5.13]{dgi} that every geodesic $\g$ representing $\zeta_1 \in H_1(L_1;\Z)$ occurs as the asymptotic orbit of at most one $\tilde{J}_{\infty}$-holomorphic small plane. It follows from \Cref{proposition:smallsplit} that \eqref{equation:smallplaneembedding} is a bijection, which means that the two moduli spaces in fact coincide. 

We therefore deduce the following proposition from \cite[Prop.\ 5.11]{dgi}:  

\prop With $\mc{M}_s(J_{\infty})$ defined as above, the following properties hold: 

\begin{itemize}
\item[(i)] The asymptotic evaluation map $\mc{M}_s(J_{\infty}) \to \G_{\zeta_1} \simeq S^1$ taking a plane to its asymptotic orbit is a diffeomorphism. 
\item[(ii)] The evaluation map $\mc{M}_s(J_{\infty}) \tms \C \to \R^4 - L$ is a smooth embedding which in fact maps into $\R^4-L=\R^4- \cup_{i=1}^n L_i$.
\end{itemize}
\eprop

The smoothing procedure from \cite[Sec.\ 5.3]{dgi} applied to the moduli space of small planes can now be used to construct a smoothly embedded solid torus $\mathcal{T} \subset \R^4 - (\bigcup_{i>1} L_i)$ with boundary equal to $L_1$. This solid torus is moreover foliated by symplectic disks by construction.  

\qed

\pf[Sketch of the proof of \Cref{proposition:enhancement}]

Observe that the characteristic distribution $\ker \omega|_{T\mathcal{T}}$ integrates to a monodromy of the disk leaf, which is a symplectomorphism that fixes the boundary. To prove \Cref{proposition:enhancement}, we wish to correct the monodromy by a suitable isotopy of $\mc{T}$.  The structure of the argument is entirely analogous to that of \cite[Sec.\ 6]{dgi}. 

As in \cite[Sec.\ 6.2]{dgi}, the monodromy can be corrected by applying the \emph{suspension construction}. This allows one to change the monodromy of $\mc{L}$ at the cost of isotoping $\mc{L}$ by some distance which depends on the size of the Hamiltonian generating the monodromy.  

The suspension construction may a priori create self-intersections of $\mc{T}$ and/or intersections of $\mc{T}$ with the other Lagrangians (the second possibility does not occur in \cite{dgi}). In order to ensure that this does not happen, one needs to ensure that there is ``enough room" around $\mc{T}$. This can be done exactly as in \cite[Sec.\ 6.1]{dgi}. More precisely, one applies a suitable \emph{inflation}. This outcome of this inflation is essentially to shrink the Lagrangian $L$ through Lagrangians. This also has the effect of shrinking the size of the Hamiltonian needed to perform the suspension, thus ensuring that there is indeed ``enough room". \epf

\section{Applications to linking}

In this section, we apply the holomorphic curve analysis of the previous section to study linking of Lagrangian tori in symplectic $4$-manifolds. In particular, we prove \Cref{theorem:lagrangian-r4-introduction} and \Cref{theorem:lagrangian-rational-introduction} from the introduction.

\subsection{Linking and disk counts}

We begin by reviewing a criterion from \cite{cote}, formulated in terms of enumerative disk counts, under which a collection of monotone Lagrangian tori satisfy an algebro-topological unlinking property. 

\defi
	Let $L$ be a monotone Lagrangian torus in a symplectic $4$-manifold $(M,\o)$. Given a class $\beta \in \pi_2(M, L)$ of Maslov number $2$, let $n(L, \beta)$ be the mod-$2$ count of $J$-holomorphic disks passing through a generic point of $L$, for a generic $J$. It can be shown by standard arguments (see \cite[Sec.\ 3.1]{auroux}) that count is well-defined, i.e. independent of $J$.
\edefi

Let us consider a collection $L_1,\dots,L_n \sub (M, \o)$ of disjoint, monotone Lagrangian tori. Let $\kappa_i$ be the monotonicity constant of $L_i$ (i.e. $\kappa_i \mu = \o \in \pi_2(M, L_i)$, where $\mu$ is the Maslov class). If $\kappa_i \leq \kappa_j$ for $i \leq j$, then it can be shown by a straightforward neck-stretching argument (see \cite[Prop.\ 4.3]{cote}) that there exists an almost-complex structure $J$ with the following property: if $n(L_i, \beta) \neq 0$, then there exists a $J$-holomorphic disk $u_{\beta}: (D^2, \d D^2) \to (M, L)$ representing the class $\beta$ such that $u_{\beta} \cap L_j= \emptyset$ for all $j \geq i$.

We easily deduce the following proposition.

\prop \label{proposition:enumerative}
Suppose that there exists a collection of classes $\mc{C}_i = \{ \beta_1,\dots,\beta_{n_i} \} \sub \pi_2(M, L_i)$ such that $n(L_i, \beta_k) \neq 0$ for all $1 \leq k \leq n_i$, and such that the image of $\mc{C}_i$ under the composition 
\eq  \label{equation:composition} \pi_2(M, L_i) \to H_2(M, L_i) \to H_1(L_i) = \pi_1(L_i) \eeq generates $\pi_1(L_i)$ as an abelian group. Then
\eq \pi_1(L_i) \to \pi_1( M- \bigcup_{j> i} L_j) \eeq is the trivial map. 
\eprop
\qed

\subsection{Unlinking Lagrangian tori in $\R^4$}  \label{subsection:unlinking-r4}

We remind the reader of the following definition. 

\defi
Let $L_1,\dots,L_n$ be disjoint, contractible Lagrangians in a symplectic manifold $(X, \o)$. We say that the $L_i$ are \emph{smoothly unlinked} if there exists a collection of disjoint embedded closed balls $B_1,\dots,B_n \sub X$ and $1$-parameter family of smooth embeddings $\Phi: [0,1] \tms (\bigsqcup_{i=1}^n L_i ) \to X$ such that, for $i=1,\dots,n$, the composition $$L_i \xhookrightarrow{L_i \mapsto \{t\} \tms L_i} [0,1] \tms \lt(\bigsqcup_{i=1}^n L_i\rt) \xrightarrow{\Phi} X$$ maps $L_i$ to $L_i \sub X$ diffeomorphically if $t=0$ and maps $L_i$ into $B_i$ if $t=1$.\footnote{We are slightly abusing notation by viewing $L_i$ both as an abstract manifold and as an embedded submanifold of $X$.}

We say that the $L_i$ are \emph{Lagrangian unlinked} (resp. \emph{Hamiltonian unlinked}) if one can choose $\Phi$ so that $\Phi(t,-)$ is a Lagrangian embedding for all $t \in [0,1]$ (resp. if $\Phi$ is induced by a global Hamiltonian isotopy of $(X,\o)$). 
\edefi 

The following theorem was stated as \Cref{theorem:lagrangian-r4-introduction} in the introduction and can be thought of as a 1-parametric h-principle for Lagrangian embeddings of a collection of tori into $\R^4$. Indeed, it implies that a collection of Lagrangian tori in $\R^4$ are Lagrangian unlinked if and only if the obvious algebro-topological obstructions vanish. This is not a priori obvious, since there could be obstructions coming from differential or symplectic topology. As a corollary of this theorem, one obtains obstructions to certain ``linked" Lagrangian embeddings of tori into $(\R^4, \o)$, as was explained in \Cref{subsection:obstructions-linked-r4}. 

\thm \label{theorem:maintheorem} 
	Let $L_1,\dots,L_n$ be disjoint Lagrangian tori in $(\R^4,\o)$.  If the inclusion $$\iota_i: L_i \to \R^4- \bigcup_{j \neq i} L_j$$ induces the trivial map on fundamental groups for $i=1,2\dots,n$, then the $L_i$ are Lagrangian unlinked. If the $L_i$ are monotone and $\kappa_i \leq \kappa_j$ for $i \leq j$, where $\kappa_i$ is the monotonicity constant of $L_i$, then it is enough to assume that $\tilde{\iota}_i: L_i \to \R^4- \bigcup_{j >i} L_j$ induces the trivial map on fundamental groups.
\ethm


We begin by arguing that it is enough to prove \Cref{theorem:maintheorem} in the non-monotone case. 

\lem \label{lemma:monotone-reduction}
Under the hypotheses of \Cref{theorem:maintheorem}, if the $L_i$ are monotone with $\kappa_i \leq \kappa_j$ and $\tilde{\iota}_i$ induces the trivial map on fundamental groups, then $\iota_i$ also induces the trivial map on fundamental groups for all $i=1,\dots,n$. 
\elem
\pf
According to \Cref{proposition:solidbound}, we can assume that $L_1$ bounds a smoothly embedded solid torus $S_1 \sub \R^4 - \bigcup_{j>1} L_j$. Let $\g_1$ be the core of this solid torus. Since the inclusion $L_1 \to \R^4- \bigcup_{j>1} L_j$ induces the trivial map on fundamental groups, it follows that $\g_1$ can be isotoped into some ball $B_1$ which is far away from the other Lagrangians. Let $\phi_1$ be this isotopy. By the isotopy extension theorem, $\phi_1$ can be extended to a global smooth isotopy which is constant near the other Lagrangians.  

Observe that $\phi_1$ also takes some small neighborhood $U_1$ of $\g_1$ into $B_1$. By contracting the leaves of the solid torus $S_1$, we can smoothly isotope $L_1$ into $U_1$ without intersecting the other tori. By concatenating this isotopy with $\phi_1$, we obtain a smooth isotopy taking $L_1$ into $B_1$ which is constant near the other Lagrangians. 

We now argue by induction. Let us suppose for $1\leq k<n$ that there exists a smooth isotopy taking $L_1,\dots, L_k$ into balls $B_1,\dots,B_k$ which are pairwise disjoint and disjoint from the other Lagrangians. We can freely assume that $B_1,\dots,B_k$ are far away from $L_{k+1},\dots,L_n$. Hence we can repeat the argument of the previous paragraph to find an isotopy taking $L_{k+1}$ into a ball $B_{k+1}$ which is disjoint from $B_1,\dots,B_k$ and from $L_{k+2},\dots,L_n$. This proves that the $L_i$ are smoothly unlinked, so in particular $\iota_i$ induces the trivial map on fundamental groups for all $i=1,\dots,n$. 
\epf

\pf[Proof of \Cref{theorem:maintheorem}] According to \Cref{proposition:enhancement}, we can assume, up to relabelling the tori and isotoping $L=L_1\cup\dots\cup L_n$ through Lagrangian tori, that $L_1$ bounds a smoothly embedded solid torus $S_1 \sub \R^4 - \bigcup_{j>1} L_j$. Moreover, $S_1$ is foliated by symplectic disks, and the monodromy map induced on these disks by the characteristic foliation is the identity. 

Let $\g_1$ be the core of $S_1$. Since the inclusion $L_1 \to \R^4- \bigcup_{j>1} L_j$ induces the trivial map on fundamental groups, it follows that $\g_1$ can be isotoped into some ball $B_1$ which is far away from the other Lagrangians. 

Arguing now as in \cite[Sec.\ 6.3]{dgi}, $L_1$ can be isotoped through Lagrangian tori to a circle bundle over $\g_1$ by contracting the leaves of the bounding solid torus. The isotopy taking $\g$ into $B_1$ can also be extended to an isotopy of this circle bundle which therefore will not intersect the other Lagrangians if the circle fibers are small enough. This produces a Lagrangian isotopy which takes $L_1$ into $B_1$ without intersecting the other Lagrangians. By the isotopy extension theorem, it can be extended to a global smooth isotopy which is constant near the other Lagrangians. 

The remainder of the proof goes by induction, exactly as in the last paragraph of the proof of \Cref{lemma:monotone-reduction}. 
\epf

As mentioned in the introduction, \Cref{theorem:maintheorem} also gives the following modest improvement of a result of the first author, who proved in \cite[Thm.\ A]{cote} that any collection of disjoint Clifford tori in $\R^4$ are \emph{smoothly} unlinked.  

\cor \label{corollary:local-unlinking} Suppose that $L_1,\dots,L_n$ are monotone Clifford tori. Then the $L_i$ are Lagrangian unlinked. \ecor
\pf 
We can assume that the $L_i$ are indexed so that $\kappa_i \leq \kappa_j$ if $i \leq j$, where $\kappa_i$ is the monotonicity constant of $L_i$. Since the $L_i$ are Clifford tori, we can consider the pair of classes $\beta_1= [D^2 \tms \{*\}] $ and $\beta_2 = [\{* \} \tms D^2]$ in $\pi_2(\R^4, L_i)$. It is well-known (see \cite[Ex.\ 4.2]{cote}) that $n(L_i, \b_1)= n(L_i, \b_2) =1$. Moreover, it is clear that the image of $\{ \beta_1, \beta_2 \}$ under \eqref{equation:composition} generates $\pi_1(L_i)$. It follows by \Cref{proposition:enumerative} that $\tilde{\iota}_i: L_i \to (\R^4- \bigcup_{j> i} L_j)$ induces the trivial map on fundamental groups. The corollary now follows from \Cref{theorem:maintheorem}. 
\epf

\rmk As explained in \cite[Ex.\ 4.10]{cote}, it is easy to construct examples of monotone Chekanov tori in $\R^4$ which are linked. \ermk

\subsection{Unlinking Lagrangian tori in symplectic rational surfaces} \label{subsection:unlinking-rational}

Let $(X, \o)$ denote either $(\bb{CP}^2, \o)$ or $(S^2 \tms S^2, a \o \oplus b \o)$, where $a, b>0$ and $\o$ always denotes the Fubini-Study form of unit volume on $\bb{CP}^2$ and $S^2= \bb{CP}^1$ respectively. If $X= \bb{CP}^2$, let $D_{\infty} \sub X$ be a fixed complex line (which we think of as the line at infinity). If $X= S^2 \tms S^2$, let $D_{\infty} \sub X$ be the divisor $S^2 \tms \{ \infty \} \cup \{ \infty \} \tms S^2$, where $\infty \in S^2$ is a fixed point.

\lem \label{lemma:displace-divisor}
	Let $L_1,\dots,L_n \sub (X, \o)$ be disjoint Lagrangian tori. There exists a global Hamiltonian isotopy of $(X,\o)$ which moves all of the $L_i$ to the complement of $D_{\infty}$. 
\elem

\pf
	The case $n=1$ and $a=b$ was proved in \cite[Sec.\ 3.3.1 and 3.3.2]{dgi}. For $n >1$, one simply sets $L := L_1 \cup \dots \cup L_n$ and the same argument goes through without change. If one additionally allows $a \neq b$, then the same argument still works; however, one must in addition argue that a generic neck-stretching sequence $\{J_n\}_{n \geq 0}$ around $L$ admits a $J_n$-holomorphic foliation in both classes $[S^2 \tms *]$ and $[* \tms S^2]$ for all $n \geq 0$. This can be proved by standard arguments; cf.\ \cite[Ex.\ 6.3(a)]{wendlcurves}. \epf


\prop \label{proposition:torus-bounding}
	Let $L_1,\dots,L_n \sub (X, \o)$ be disjoint Lagrangian tori. Then, after possibly isotoping the $L_i$ through disjoint Lagrangian tori, one of the $L_i$ bounds a solid torus foliated by symplectic disks such that the monodromy map induced by the characteristic foliation is the identity. 
\eprop

\pf According to \Cref{lemma:displace-divisor}, we can assume that $L= L_1 \cup\dots\cup L_n \sub X- D_{\infty}$.  For $X=S^2 \tms S^2$ (resp.\ for $X= \bb{CP}^2$) observe that $X-D_{\infty}$ is a polydisk (resp.\ a ball) which embeds symplectically into $(\R^4, \o)$. We can therefore view $L$ as a union of Lagrangian tori in $(\R^4, \o)$. \Cref{proposition:torus-bounding} thus follows from \Cref{proposition:enhancement} (combined with a suitable rescaling which ensures that the isotopy of tori furnished by \Cref{proposition:enhancement} stays inside $X-D_{\infty} \sub \R^4$).  \epf

\thm \label{theorem:minimal-unlinked}
	Let $L_1,\dots,L_n \sub (X, \o)$ be disjoint Lagrangian tori. If the inclusion $$\iota_i: L_i \to X- \cup_{j \neq i} L_j$$ induces the trivial map on fundamental groups for $i=1,2\dots,n$, then the $L_i$ are Lagrangian unlinked. 
\ethm

\pf 
	Applying \Cref{proposition:torus-bounding}, the proof is identical to that of \Cref{theorem:maintheorem}. 
\epf

We move on to discussing linking of tori in general rational symplectic $4$-manifolds. Let us begin by reviewing the notion of symplectic blowup and blowdown, closely following \cite[Chap.\ 7]{mcduff-sal-intro}. 

Let $\mc{O}(1)$ be the tautological line bundle over $\bb{CP}^{n-1}$. It admits two natural projections $\op{pr}: \mc{O}(1) \to \bb{CP}^{n-1}$ and $\pi: \mc{O}(1) \to \C^n$, where $\op{pr}$ is the bundle projection onto the base and $\pi(\ell, x)=x$ for $x \in \C^n$ and $x \in \ell \in \bb{CP}^{n-1}$. For $\l >0$, let $\tilde{\o}_{\l}:= \pi^* \o_0 + \l^2 \op{pr}^* \o_{\op{FS}}$, where $\o_0$ is the standard symplectic form on $\C^n$ and $\o_{\op{FS}}$ is the Fubini-Study form on $\bb{CP}^{n-1}$, normalized so that $\bb{CP}^1$ has area $\pi$.  

For $r>0$, let $B(r) \sub \C^n$ be the (closed) ball of radius $r$ centered at the origin. Let $L(r)= \pi^{-1}(B(r))$.  Letting $Z \sub \mc{O}(1)$ be the zero section and fixing the symplectic form $\tilde{\o}_{\l}$ on $\mc{O}(1)$, it can be shown that $L(r)- Z$ is symplectomorphic to $B(\sqrt{\l^2+ r^2}) - B(r)$. 

Let us now consider a symplectic $2n$-manifold $(X, \o)$ and a sphere $\S \sub X$ whose normal bundle is isomorphic to $\mc{O}(1)$. (If $2n=4$, this is equivalent to the condition that $\S$ has self-intersection $-1$.) It follows from the symplectic neighborhood theorem that $\S$ has a neighborhood which is symplectomorphic to $L(\e)$ for some $\e>0$. The \emph{blowdown} of $X$ along $\S$ is now constructed by replacing $\S$ by $B(\e)$, using the fact that $L(\e)- Z$ is symplectomorphic to $B(\sqrt{\l^2+ \e^2}) - B(\e)$. The resulting manifold is uniquely defined up to symplectomorphism.  

The inverse operation is called the \emph{blowup} and is defined as follows. Given a symplectic embedding $B(\l) \hookrightarrow (X, \o)$, it can always be extended to a symplectic embedding $B(\sqrt{\l^2+ \e^2})$ for some $\e>0$. One now replaces $B(\sqrt{\l^2+ \e^2})$ with $L(\e)$. This operation is also well-defined up to symplectomorphism, and in fact depends only on the symplectic isotopy class of the embedding $B(\l) \hookrightarrow (X, \o)$.

\defi
	A closed symplectic $4$-manifold $(M,\o)$ is said to be \emph{rational} if it can be obtained from $(\bb{CP}^2, \l \o)$, $\l>0$, by some sequence of blowups and blowdowns. 
\edefi

In addition to $\bb{CP}^2$ and $S^2 \tms S^2$, well-studied examples of rational symplectic $4$-manifolds include the monotone symplectic del Pezzo surfaces. It follows from classification results of Gromov and McDuff \cite[Thm.\ 2.1]{lalonde-mcduff} that given a complex projective surface $(X, j)$ which is rational in the algebro-geometric sense (i.e. it contains a rational curve of non-negative self-intersection), any choice of $j$-compatible symplectic form $\o$ makes $(X, \o)$ into a rational symplectic $4$-manifold.  In particular, given a projective embedding $X \hookrightarrow \bb{CP}^N$, we find that $(X, \o_{\op{FS}}|_X)$ is a rational symplectic $4$-manifold.  

\prop[Gromov, McDuff] \label{proposition:blowdown}
	Suppose that $(M,\o)$ is rational. There exists a disjoint collection $E_1,\dots,E_l$ such that if we let $(\ov{M}, \ov{\o})$ denote the blowdown along the $E_i$, then $(\ov{M}, \ov{\o})$ is symplectomorphic to either $\bb{CP}^2$ (with the standard symplectic form) or to $(S^2 \tms S^2, a \o \oplus b \o)$ for $a,b >0$. 
\eprop

\pf
	McDuff proved in \cite[Thm.\ 1.1]{mcduff-jams} that there exists a collection $E_1,\dots,E_l$ of disjoint exceptional spheres such that $(\ov{M}, \ov{\o})$ is minimal. She also also showed in \cite[Thm.\ 1.2]{mcduff-jams} that the class of symplectic manifolds which contain a symplectically embedded sphere of non-negative self-intersection is closed under blowups and blowdowns. So in particular, $(\ov{M}, \ov{\o})$ admits such a sphere. 
	The minimal symplectic manifolds which admit a symplectically embedded sphere of non-negative self-intersection are fully classified by combined work of Gromov and McDuff; see for instance \cite[Thm.\ D]{wendlcurves}. Since $(\ov{M}, \ov{\o})$ is simply-connected, the only possibilities are those stated in the proposition. \epf

It will be convenient to record the following definition. 

\defi \label{definition:exceptional-sphere}
	A symplectically embedded sphere $u: S^2 \to (X^4, \o)$ of self-intersection $-1$ is said to be an \emph{exceptional sphere}. Note that $c_1(u):= \langle c_1(TX), u_*[S^2] \rangle= c_1(N_u)+ c_1(T_u)= [u] \cdot [u] + \chi(S^2)$, which implies that $c_1(u)=1$. 
\edefi

We now have the following proposition, whose proof largely follows \cite{evans}. 

\prop \label{proposition:divisor-push}
	Let $L_1,\dots,L_n$ be Lagrangian tori in a symplectic $4$-manifold $(X, \o)$. Given a collection $E_1,\dots,E_l$ of disjoint exceptional spheres, there exists a Hamiltonian isotopy of $(X,\o)$ which takes the $L_i$ into the complement of $E_1 \cup\dots\cup E_l$.  
\eprop

\pf 
	Choose a compatible complex structure $J_0$ so that the $E_i$ are $J_0$-holomorphic. After possibly perturbing $J_0$, we can assume that it is generic: indeed, by \Cref{definition:exceptional-sphere} and automatic transversality (see \cite[Lem.\ 3.3.3]{mcduff-sal-jcurves}), the exceptional curves are transversely cut-out, so they survive under small perturbations of the complex structure. It also follows from positivity of intersection that the curves remain disjoint when perturbing the complex structure. 

Let $L= \bigcup_{i=1}^n L_i$ and let $\tilde{\eta}= (g, \phi, J_{\op{cyl}}, J_{\op{std}}, J_{\infty}) \in \mc{D}(X, \o; L)$ be a neck-stretching datum satisfying the following properties (cf.\ \Cref{subsection:neckstretch}):

\begin{itemize}
\item $g$ is a rescaling of the Euclidean metric on $\bb{T}^2= \R^2/\Z^2$.
\item $\phi=(\phi_1,\dots,\phi_n)$ is a Weinstein embedding, where $\phi_i$ carries $N(0_{\bb{T}^2})$ into $N(L_i)$. Here $N(0_{\bb{T}^2}) \sub T^* \bb{T}^2$ is a neighborhood of the zero section and $N(L_i) \sub X$ is a neighborhood of $L_i$, where $N(L_i) \cap N(L_j)= \emptyset$ if $i \neq j$.
\item $S^*_{\leq 4,g} \bb{T}^2 \sub N(0_{\bb{T}^2})$.
\item $J_{\op{std}}$ and $J_{\op{cyl}}$ are defined as in \cite[Sec.\ 4]{dgi}.
\item The pullback of $J_{\infty}$ under $\phi_i$ coincides with $J_{\op{cyl}}$ on $S^*_{\leq 4,g} - 0_{\bb{T}^2}$. 
\item $J_{\infty}$ is generic in the complement of the $N(L_i)$.
\end{itemize}

Now stretch the neck along $L= \bigcup_{i=1}^n L_i$ using the datum $\tilde{\eta}$. This produces a sequence of complex structures $\{J_n\}_{n=0}^{\infty}$ which are independent of $n$ outside of a small neighborhood of $L$. Note that $J_{n}$ is regular for simply-covered genus $0$ curves due to our genericity assumption on $J_{\infty}$. For each $n \geq 0$, choose a generic path $\{J_s\}_{s \in [n, n+1]}$ joining $J_n$ to $J_{n+1}$. According to \cite[Thm.\ 5.1]{wendlcurves}, the $E_i$ are represented by unique, disjoint $J_s$-holomorphic curves which vary smoothly with $s \in [n, n+1]$. Let $E_i^s$ be the unique $J_s$-holomorphic sphere isotopic to $E_i$. 

We claim that the $E_i^n$ are disjoint from $L$ for $n$ large enough. Suppose for contradiction that this is not the case. After possibly relabelling the $E_i$, we can assume that $E_1^n$ intersects $L$ for all but finitely many $n \in \N$. By the SFT compactness, the $E_1^n$ converge to a holomorphic building $u_{\infty}$ which must have a non-empty bottom level. 

A straightforward topological argument shows that the building $u_{\infty}$ must have two planes $u_1, u_2$, which must live in the top level. According to the index formula in \Cref{proposition:compactformula}, we have $\op{ind}(u_i)= -1 + 2c_1^{\Phi}(u_i)$. In particular, we must have $c_1^{\Phi}(u_i) > 0$ since the $u_i$ are transversely cut out (if the $u_i$ are multiple covered, we can pass to the cover and run the same argument). On the other hand, the relative first Chern number is additive under gluing of building components. Moreover, as noted in \cite[Lem.\ 3.1]{dgi}, $c_1^{\Phi}(v)=0$ for any curve $v$ in a middle or bottom level. It follows that $c_1(E_i) = c_1(u_{\infty}) \geq c_1^{\Phi}(u_1) + c_1^{\Phi}(u_2) \geq 2$. This is a contradiction since $c_1(E_i)= 1$, as observed in \Cref{definition:exceptional-sphere}.

We have produced an isotopy of the $E_i$ through symplectic spheres which becomes disjoint from $L$. It is well-known (see for instance \cite[Prop.\ 0.3]{siebert-tian}) that such an isotopy extends to a global Hamiltonian isotopy $\{\phi_t\}_{t\in [0,1]}$. Observe finally that $\phi^{-1}_1(L)$ is disjoint from the $E_i$, which completes the proof. \epf

\rmk 
	Strictly speaking, we need the $J_n$ to be generic in order to appeal to \cite[Thm.\ 5.1]{wendlcurves} in the above argument. However, the proof of \cite[Thm.\ 5.1]{wendlcurves} only needs the fact that $J_n$ is regular for simply-covered curves of genus $0$. \ermk

\lem \label{lemma:blowup-complement}
	 Let $(\tilde{X}, \tilde{\o})$ be obtained by blowing up $(X, \o)$ along disjoint balls, where $(X, \o)$ denotes either $(\bb{CP}^2, \o)$ or $(S^2 \tms S^2, a \o \oplus b \o)$ for $a, b>0$. Given disjoint Lagrangian tori $L_1,\dots,L_n \sub (\tilde{X}, \tilde{\o})$, one of them bounds a solid torus which is disjoint from the other tori and from the exceptional divisors.
\elem

\pf
	According to \Cref{proposition:divisor-push}, we can assume after moving the tori by a global Hamiltonian isotopy that they do not intersect the exceptional spheres of $(\tilde{X}, \tilde{\o})$. They therefore correspond to tori in $(X, \o)$ via blowdown. According to \Cref{proposition:torus-bounding}, one of these tori, say $L_1$, bounds a solid torus $S_1$. Observe that the blowup is a symplectomorphism away from a collection of disjoint balls $B_1, \dots, B_l \sub X$, and that the $L_i$ are contained in $X- \bigcup_{j=1}^l B_j$. Since the $B_j$ are contractible, we may assume after possibly isotoping $S_1$ while keeping its boundary fixed that $S_1$ is contained in  $X- \bigcup_{j=1}^l B_j - \bigcup_{i=2}^n L_i$. Hence $S \cup L_2 \cup \dots L_n$ lifts to $(\tilde{X}, \tilde{\o})$. \epf

\pf[Proof of \Cref{proposition:all-tori-isotopic}] Let $L_1, L_2$ be Lagrangian tori in $(M, \o)$. Applying \Cref{lemma:blowup-complement} in the case $n=1$ to each $L_i$, we find that $L_i$ bounds a solid tori $S_i \sub M$. Let $B \sub M$ be a fixed Darboux ball. Since blowups and blowdown do not affect the fundamental group, we have that $\pi_1(M)=0$. Hence each $S_i$ can be isotoped into $B$. The claim now follows from the straightforward observation that all embeddings of a (closed) solid torus into a $4$-dimensional ball are isotopic.  
\epf

\pf[Proof of \Cref{corollary:vanishing-linking-class-rational}] 
We claim first that all Lagrangian tori in $\R^4, \bb{CP}^2$ and $(S^2 \tms S^2, a \o \oplus b \o)$ are Lagrangian isotopic. For a torus $L  \sub (\R^4, \o)$, this is \cite[Thm.\ A]{dgi}. The general case can be reduced to $(\R^4, \o)$ by applying \Cref{lemma:displace-divisor} to move $L$ into a ball or polydisk. 

Observe next that the linking class is preserved under Lagrangian isotopy. Since a product torus in a Darboux ball clearly has vanishing linking class, we conclude that all Lagrangian tori in $\R^4, \bb{CP}^2$ and $(S^2 \tms S^2, a \o \oplus b \o)$ have vanishing linking class.

We now consider an arbitrary rational symplectic $4$-manifold $(M, \o)$ and let $L \sub M$ be a Lagrangian torus. By \Cref{proposition:blowdown}, $(M, \o)$ is obtained by blowing up a symplectic manifold $(\ov{M}, \ov{\o})$ along disjoint balls $B_1,\dots,B_l$, where $(\ov{M}, \ov{\o})$ is symplectomorphic to $\bb{CP}^2$ or $(S^2 \tms S^2, a \o \oplus b \o)$. 

By the definition of the blowdown, there is a symplectomorphism $\pi: M - (\bigcup_{i=1}^n Z_i) \to \ov{M}- (\bigcup_{i=1}^n B_i)$. \Cref{lemma:blowup-complement} in the case $n=1$ shows that $L$ bounds a smoothly embedded solid torus $S$ which is disjoint from the exceptional divisors $Z_1\cup\dots\cup Z_n$. This implies that $\pi: L \to \pi(L)$ preserves the nullhomologous framing on $L$. Since $\pi$ is a symplectomorphism, it automatically preserves the Lagrangian framing. Hence $\pi: L \to \pi(L)$ preserves the linking class. Since we have already established that $\pi(L)$ has vanishing linking class, it follows that $L$ does too. 
\epf

\pf[Proof of \Cref{theorem:lagrangian-rational-introduction}]
If $(M, \o)$ is minimal, and hence is symplectomorphic to $(\bb{CP}^2, \o)$ or $(S^2 \tms S^2, a\o \oplus b \o)$, then the result follows from \Cref{theorem:minimal-unlinked}. 

In general, it follows from \Cref{proposition:blowdown} that $(M, \o)$ is a blowup of $(\bb{CP}^2, \o)$ or $(S^2 \tms S^2, a \o \oplus b \o)$ along a union of disjoint balls. By \Cref{lemma:blowup-complement}, we may assume after possibly relabelling the tori that $L_1$ bounds a smoothly embedded solid torus which is contained in the complement of $L_2 \cup \dots \cup L_n$. The remainder of the proof is now analogous to that of \Cref{theorem:maintheorem}: since $\iota_1$ induces the trivial map on fundamental groups, we can isotope $L_1$ into a small ball $B_1$ which is disjoint from the other tori. We now proceed by induction, repeating the same argument with one fewer torus until the process terminates. \epf

\end{document}